\documentclass[11pt]{article} 
\usepackage[margin=1.4in]{geometry}
\usepackage[utf8]{inputenc}
\usepackage[T1]{fontenc}
\usepackage[english]{babel}
\usepackage{lmodern} 
\usepackage{graphicx}
\usepackage{wrapfig}
\usepackage{subfig}
\usepackage{caption}
\usepackage{paralist}
\usepackage{nicefrac}

\graphicspath{{./figures/}}

\newcommand{\titel}{On the Circumference of Essentially 4-connected Planar Graphs}

\usepackage[usenames]{xcolor} 
\definecolor{hellblau}{rgb}{0.2,0.4,1} 
\definecolor{dunkelblau}{rgb}{0,0,0.8}
\definecolor{dunkelgruen}{rgb}{0,0.5,0}
\usepackage[
	pdftex,
	colorlinks,
	linkcolor=dunkelblau,
	urlcolor=dunkelblau,
	citecolor=dunkelgruen,
	bookmarks=true,
	linktocpage=true,
	pdfsubject={}
]{hyperref} 
\usepackage{pdfpages}
\usepackage{tikz}

\usepackage{amsmath} 
\usepackage{amsthm} 
\usepackage{amsfonts} 
\theoremstyle{plain} 
	\newtheorem{satz}{Satz}[] 
	\newtheorem{theorem}[satz]{Theorem}
	\newtheorem{lemma}[satz]{Lemma}

\theoremstyle{remark} 
\theoremstyle{definition} 

\let\circ\relax
\DeclareMathOperator\circ{circ}

\usepackage[normalem]{ulem}

\begin{document}
	\title{\titel}
	\author{}
		\author{Igor Fabrici\thanks{Partially supported by DAAD, Germany (as part of BMBF) and by the Ministry of Education, Science, Research and Sport of the Slovak Republic within the project 57320575.} \thanks{Partially supported by Science and Technology Assistance Agency under the contract No.\ APVV-15-0116 and by the Slovak VEGA Grant 1/0368/16.} \thanks{Institute of Mathematics, Pavol Jozef Šafárik University, Košice}
		\and Jochen Harant\footnotemark[1] \thanks{Institute of Mathematics, TU Ilmenau}
		\and Samuel Mohr\footnotemark[1] \footnotemark[4] \thanks{Gefördert durch die Deutsche Forschungsgemeinschaft (DFG) -- 327533333 und 270450205; partially supported by the grants 327533333 and SCHM 3186/1\mbox{-}1 (270450205) from the Deutsche Forschungsgemeinschaft (DFG, German Research Foundation).}
		\and Jens M.\ Schmidt\footnotemark[1] \footnotemark[4] \footnotemark[5]}
	\date{}
	\maketitle

\begin{abstract}
A planar graph is \emph{essentially $4$-connected} if it is 3-connected and every of its 3-separators is the neighborhood of a single vertex. Jackson and Wormald proved that every essentially 4-connected planar graph $G$ on $n$ vertices contains a cycle of length at least $\frac{2n+4}{5}$, and this result has recently been improved multiple times.

In this paper, we prove that every essentially 4-connected planar graph $G$ on $n$ vertices contains a cycle of length at least $\frac{5}{8}(n+2)$. This improves the previously best-known lower bound $\frac{3}{5}(n+2)$.
\end{abstract}

\section{Introduction}
The \emph{circumference} $\circ(G)$ of a graph $G$ is the length of a longest cycle of $G$. Originally being the subject of Hamiltonicity studies, essentially 4-connected planar graphs and their circumference have been thoroughly investigated throughout literature. Jackson and Wormald~\cite{Jackson1992} proved that $\circ(G)\geq \frac{2n+4}{5}$ for every essentially 4-connected planar graph $G$ on $n$ vertices. An upper bound is given by an infinite family of essentially 4-connected planar graphs $G$ such that $\circ(G)=\frac{2}{3}(n+4)$~\cite{Fabrici2016}. Fabrici, Harant and Jendroľ~\cite{Fabrici2016} improved recently the lower bound to $\circ(G)\ge \frac{1}{2}(n+4)$; this result in turn was strengthened to $\circ(G)\ge \frac{3}{5}(n+2)$ in~\cite{Fabrici2018}. It remained an open problem whether every essentially $4$-connected planar graph $G$ on $n$ vertices satisfies $\circ(G) > \frac{3}{5}(n+2)$.

In this paper, we present the following result.

\begin{theorem}\label{thm:main}
Every essentially 4-connected planar graph $G$ on $n$ vertices contains a cycle of length at least $\frac{5}{8}(n+2)$. If $n \geq 16$, $\circ(G) \geq \frac{5}{8}(n+4)$.
\end{theorem}

This result encompasses most of the results known for the circumference of essentially 4-connected planar graphs (some of which can be found in~\cite{Fabrici2016,Grunbaum1976,Zhang1987}). In particular, it improves the bound $\circ(G)\ge \frac{13}{21}(n+4)$ that has been given in~\cite{Fabrici2016} for the special case that $G$ is maximal planar for sufficiently large $n$ (in fact, for every $n \geq 16$, as explained in Section~\ref{sec:remarks}).

\section{Preliminaries}
Throughout this paper, all graphs are simple, undirected and finite.
For a vertex $x$ of a graph $G$, denote by $\deg_G(x)$ the degree of $x$ in $G$. For a vertex subset $A \subseteq V$, let the \emph{neighborhood} $N_G(A)$ of $A$ consist of all vertices in $V-A$ that are adjacent to some vertex of $A$. For vertices $v_1,v_2,\dots,v_i$ of a graph $G$, let $(v_1,v_2,\dots,v_i)$ be the path of $G$ that visits the vertices in the given order. We omit subscripts if the graph $G$ is clear from the context.

A \emph{separator $S$} of a graph $G$ is a subset of $V$ such that $G-S$ is disconnected; $S$ is a $k$-\emph{separator} if $|S|=k$. 
A separator $S$ is \emph{trivial} if at least one component of $G-S$ is a single vertex, and \emph{non-trivial} otherwise. Let a graph $G$ be \emph{essentially 4-connected} if $G$ is $3$-connected and every 3-separator of $G$ is trivial. It is well-known that, for every 3-separator $S$ of a 3-connected planar graph $G$, $G - S$ has exactly two components.

A cycle $C$ of a graph $G$ is \emph{isolating} (sometimes also called \emph{outer-independent}) if every component of $G-V(C)$ is a single vertex that has degree three in $G$. An edge $xy$ of a cycle $C$ of $G$ is \emph{extendable} if $x$ and $y$ have a common neighbor in $G-V(C)$. For example, Figure~\ref{fig:2face2} depicts (a part) of an isolating cycle $C$ for which the edge $yz$ becomes extendable after contracting the edge $zu$. According to Whitney~\cite{Whitney1932}, every 3-connected planar graph has a unique embedding into the plane (up to flipping). Hence, we assume in the following that the embeddings of such graphs are fixed.

\section{Proof of Theorem~\ref{thm:main}}
Let $G$ be an essentially 4-connected plane graph. It is well-known that every 3-connected plane graph on at most 10 vertices is Hamiltonian~\cite{Dillencourt1996}; thus, for $4\leq n\leq 10$, this implies $\circ(G) = n \geq \frac{5}{8}(n+2)$. 
Since these graphs contain in particular the essentially 4-connected plane graphs on at most 10 vertices, we assume $n \geq 11$ from now on. 
For $n \geq 11$, it was shown in~\cite[Lemma~4(ii)]{Fabrici2016} that $G$ contains an isolating cycle of length at least $8$. Let $C$ be a longest such isolating cycle of length $c := |E(C)| \geq 8$. We will show that $c \geq \frac{5}{8}(n+2)$, so that $C$ is a cycle of the desired length.

Clearly, $C$ contains no extendable edge $xy$, as otherwise one could find a longer such cycle by replacing $xy$ in $C$ with the path $(x,v,y)$, where $v \notin V(C)$ is a common neighbor of $x$ and $y$.
Let $V^-$ be the subset of vertices of $V$ that are contained in the open set of $\mathbb{R}^2-C$ that is bounded (hence, strictly inside $C$), and let $V^+ := V-V(C)-V^-$. We assume that $|V^-| \geq 1 \leq |V^+|$, since otherwise we are done, as then $c \geq \frac{2}{3}(n+2)$ is implied by~\cite[Lemma~5]{Fabrici2016}. Let $H$ be the plane graph obtained from $G$ by deleting all chords of $C$ (i.\,e., all edges $xy \in E-E(C)$ satisfying $x,y \in V(C)$) and let $H^- := H-V^+$ and $H^+ := H-V^-$. A face of $H$ is called \emph{minor} if it is incident to exactly one vertex of $V^- \cup V^+$, and \emph{major} otherwise. Let $M^-$ and $M^+$ be the sets of minor faces in $H^-$ and $H^+$, respectively. For example, in Figure~\ref{fig:2face2}, we have $a \in V^-$, $b \in V^+$, $f \in M^-$ and $f' \in M^+$.

Note that a face $f$ of $H$ is incident to no vertex of $V^- \cup V^+$ if and only if it is bounded by $C$ (i.e., if $f$ is either the region inside or outside $C$). Since we assumed $|V^-| \geq 1 \leq |V^+|$, our definition of minor faces coincides with the one of~\cite{Fabrici2018}, so that we can use the following inequality.

\begin{lemma}[{\cite{Fabrici2018}, Inequality~(i)}]\label{lem:minor}
$|M^- \cup M^+| \geq |V^- \cup V^+|+2$.
\end{lemma}

In $H$, an edge $e$ of $C$ is incident with exactly two faces $f$ and $f'$ of $H$. In this case we say $f'$ is \emph{opposite} to $f$ with respect to $e$. A face $f$ of $H$ is called $j$-\emph{face} if it is incident with exactly $j$ edges of $C$; the edges of $C$ that are incident with $f$ are called $C$-\emph{edges} of $f$. Since $C$ does not contain an extendable edge, we have $j \geq 2$ for every minor $j$-face of $H$. For two faces $f$ and $f'$ of $H$, let $m_{f,f'}$ be the number of common $C$-edges of $f$ and $f'$.

\bigskip
If we can prove 
\begin{align}\label{eq:facestovertices}
2c \geq \frac{10}{3}|M^- \cup M^+|,
\end{align}

\noindent then Theorem~\ref{thm:main} follows directly from the inequality $|M^- \cup M^+| \geq n-c+2$ of Lemma~\ref{lem:minor}.
We charge every $j$-face of $H$ with weight $j$ (and thus have a total charge of weight $2c$) and discharge these weights in $H$ by applying the following set of rules exactly once. In order to prove Inequality~\eqref{eq:facestovertices}, we will aim to prove that every minor face of $H$ has weight at least \nicefrac{10}{3} after the discharging.

\begin{description}
	\item[Rule R1:] Every major face $f$ of $H$ sends weight $m_{f,f'}$ to every minor face $f'$ opposite to $f$.
	\item[Rule R2:] Every minor face $f$ of $H$ sends weight $\frac23 m_{f,f'}$ to every minor 2-face $f'$ opposite to $f$.
	\item[Rule R3:] Every minor face $f$ of $H$ sends weight $1$ to every minor 3-face $f'$ that is opposite to $f$ with respect to the middle $C$-edge of $f'$.
	\item[Rule R4:] Let $f_1$ be a minor 4-face that has an opposite minor $j$-face $f$ satisfying $j \geq 4$ and $m_{f_1,f} = 2$, as well as an opposite minor 2- or 3-face $f_2$ satisfying $m_{f_1,f_2} = 2$. Then $f$ sends weight \nicefrac{2}{3} to $f_1$.
	\item[Rule R5:] Let $f_1$ be a minor 5-face that has an opposite minor $j$-face $f$ satisfying $j \geq 4$ and $m_{f_1,f} = 2$, as well as two opposite minor 2-faces. Then $f$ sends weight \nicefrac{1}{3} to $f_1$.
\end{description}

For example, in Figure~\ref{fig:2face2}, both faces $f$ and $f'$ would send weight $2/3$ to each other according to Rule~R2, which effectively cancels the exchange of weights. Rules~R2 and~R3 may be seen as a refinement of the two rules given in~\cite{Fabrici2018}; for that reason, some of the early cases about minor 2- and 3-faces in the following case distinction will be similar as in~\cite{Fabrici2018}.

Let $w$ denote the weight function on the set $F(H)$ of faces of $H$ after Rules~R1--R5 have been applied. Clearly, $\sum_{f\in F(H)} w(f) = 2c$ still holds. In order to prove that the weight $w(f)$ of every minor face $f$ of $H$ is at least \nicefrac{10}{3} and no major face has negative weight, we distinguish several cases. For most of them, we construct a cycle $\overline{C}$ that is obtained from $C$ by replacing a subpath of $C$ with another path. In such cases, $\overline{C}$ will be an isolating (which is easy to verify due to $V(C) \subseteq V(\overline{C})$) cycle of $G$ that is longer than $C$ (we say $C$ is \emph{extended}); this contradicts the choice of $C$ and therefore shows that the considered case cannot occur. Note that the vertices of $C$ that are depicted in the following figures are pairwise non-identical, because $c \geq 8$; in the rare figures that show more than $8$ vertices of $C$, $C$ has always at least the number of vertices shown.

\bigskip\noindent
\textbf{Let $f\in F(H)$.}%
\begin{description}
\item[Case 1:] \textit{$f$ is a major $j$-face for any $j$.}

Initially, $f$ is charged with weight $j$. By Rule~R1, $f$ sends for every of its $C$-edges weight at most $1$ to an opposite face. We conclude $w(f) \geq 0$.

\item[Case 2:] \textit{$f$ is a minor $2$-face \emph{(see Figure~\ref{fig:2face})}.}

Let $xy$ and $yz$ be the $C$-edges of $f$ and let $a$ be the vertex of $V-V(C)$ that is incident with $f$. The face $f$ is initially charged with weight $2$ and gains weight at least \nicefrac{4}{3} by~R1 and~R2. If $f$ does not send any weight to other faces, this gives $w(f) \geq \nicefrac{10}{3}$, so assume that $f$ sends weight to some face $f' \neq f$.

\begin{figure}[!ht]
\begin{center}
	\begin{tikzpicture}[scale=0.25, -, 
	vertex/.style={circle,fill=black,draw,minimum size=5pt,inner sep=0pt}]
	\node[vertex,label=270:{ $x$}] (x) at (0,0) [] {$ $};
	\node[vertex,label=270:{ $y$}] (y) at (4,0) [] {$ $};
	\node[vertex,label=270:{ $z$}] (z) at (8,0) [] {$ $};
	\node[vertex,label=90:{ $a$}] (a) at (4,7) [] {$ $};
	\draw (4,4) node {$f$};
	\node (p1) at(0,0) [] {$ $};
	\node (p2) at(8,0) [] {$ $};
	\node (q1) at(-4,1) [] {$ $};
	\node[label=0:{ $C$}] (q2) at(12,1) [] {$ $};
	\path[-,out=-30, in=180] (q1) edge (p1);
	\draw (x) -- (z);
	\path[-,out=0, in=-150] (p2) edge (q2);
	\draw (x) -- (a);
	\draw (z) -- (a);
	\end{tikzpicture}
\end{center}
\caption{Case~2}\label{fig:2face}
\end{figure}

According to~R1--R5, $f'$ is opposite to $f$ and either a minor $2$-face or a minor $3$-face of $H$. Without loss of generality, let $f'$ be opposite to $f$ with respect to the edge $yz$. We distinguish the following subcases.

\begin{description}
\item[Case 2a:] \textit{$f'$ is a minor $2$-face and $xy$ is a $C$-edge of $f'$.}

Then $\{x,z\}$ is the neighborhood of $y$ in $G$, which contradicts the $3$-connectivity of $G$.

\item[Case 2b:] \textit{$f'$ is a minor $2$-face and $xy$ is not a $C$-edge of $f'$ \emph{(see Figure~\ref{fig:2face2})}.}

Then a longer isolating cycle $\overline{C}$ is obtained from $C$ by replacing the path $(x,y,z,u)$ with the path $(x,a,z,y,b,u)$ (see Figure~\ref{fig:2face2}), which contradicts the choice of $C$.

\begin{figure}[!htb]
\begin{minipage}[b]{.5\textwidth}
\begin{center}
	\begin{tikzpicture}[scale=0.25, -, 
	vertex/.style={circle,fill=black,draw,minimum size=5pt,inner sep=0pt}]
	\node[vertex,label=270:{ $x$}] (x) at (0,0) [] {$ $};
	\node[vertex,label=270:{ $y$}] (y) at (4,0) [] {$ $};
	\node[vertex,label=270:{ $z$}] (z) at (8,0) [] {$ $};
	\node[vertex,label=270:{ $u$}] (u) at (12,0) [] {$ $};
	\node[vertex,label=90:{ $a$}] (a) at (4,7) [] {$ $};
	\draw (4,4) node {$f$};
	\draw (8,-4.5) node {$f'$};
	\node[vertex,label=-90:{ $b$}] (b) at (8,-7) [] {$ $};
	\node (p1) at(0,0) [] {$ $};
	\node (p2) at(12,0) [] {$ $};
	\node (q1) at(-4,1) [] {$ $};
	\node[label=0:{ $C$}] (q2) at (16,1) [] {$ $};
	\path[-,out=-30, in=180] (q1) edge (p1);
	\draw (x) -- (u);
	\path[-,out=0, in=-150] (p2) edge (q2);
	\draw (x) -- (a);
	\draw (z) -- (a);
	\draw (y) -- (b);
	\draw (u) -- (b);
	\end{tikzpicture}
\end{center}
\caption{Case~2b}\label{fig:2face2}
\end{minipage}\hfil %
\begin{minipage}[b]{.5\textwidth}
\begin{center}
	\begin{tikzpicture}[scale=0.25, -, 
	vertex/.style={circle,fill=black,draw,minimum size=5pt,inner sep=0pt}]
	\node[vertex,label=270:{ $x$}] (x) at (0,0) [] {$ $};
	\node[vertex,label=270:{ $y$}] (y) at (4,0) [] {$ $};
	\node[vertex,label=270:{ $z$}] (z) at (8,0) [] {$ $};
	\node[vertex,label=270:{ $u$}] (u) at (12,0) [] {$ $};
	\node[vertex,label=90:{ $a$}] (a) at (4,7) [] {$ $};
	\draw (4,4) node {$f$};
	\draw (6,-4.5) node {$f'$};
	\node[vertex,label=-90:{ $b$}] (b) at (6,-7) [] {$ $};
	\node (p1) at(0,0) [] {$ $};
	\node (p2) at(12,0) [] {$ $};
	\node (q1) at(-4,1) [] {$ $};
	\node[label=0:{ $C$}] (q2) at(16,1) [] {$ $};
	\path[-,out=-30, in=180] (q1) edge (p1);
	\draw (x) -- (u);
	\path[-,out=0, in=-150] (p2) edge (q2);
	\path[dashed,out=-45, in=-135] (y) edge (u);
	\draw (x) -- (a);
	\draw (z) -- (a);
	\draw (x) -- (b);
	\draw (u) -- (b);
	\end{tikzpicture}
\end{center}
\caption{Case~2c}\label{fig:2face3}
\end{minipage}
\end{figure}

\item[Case 2c:] \textit{$f'$ is a minor $3$-face \emph{(see Figure~\ref{fig:2face3})}.}

Since we assumed that $f$ sends weight to $f'$, one $C$-edge of $f$, say without loss of generality $yz$, is the middle $C$-edge of $f'$, according to~R3. The edge $yu$ (see Figure~\ref{fig:2face3}) exists in $G$ (but not in $H$, as $H$ does not contain chords of $C$), because otherwise $d_G(y)=2$, which contradicts that $G$ is 3-connected. Then $\overline{C}$ is obtained from $C$ by replacing the path $(x,y,z,u)$ with the path $(x,a,z,y,u)$.
\end{description}

\item[Case 3:] \textit{$f$ is a minor $3$-face \emph{(see Figure~\ref{fig:3face})}.}

Then $f$ is initially charged with weight~3 and gains weight at least~1 by~R1 and~R3. If $f$ sends weight at most \nicefrac{2}{3} to other faces, this gives $w(f) \geq \nicefrac{10}{3}$, so assume that $f$ sends weight more than~\nicefrac{2}{3}. Since all weights are multiples of \nicefrac{1}{3}, $f$ has to send weight at least~\nicefrac{3}{3}. In particular, this implies that Rule~R2 or~R3 applies on $f$.

\begin{figure}[!htb]
\begin{center}
	\begin{tikzpicture}[scale=0.25, -, 
	vertex/.style={circle,fill=black,draw,minimum size=5pt,inner sep=0pt}]
	\node[vertex,label=270:{ $v$}] (v) at (-4,0) [] {$ $};
	\node[vertex,label=270:{ $x$}] (x) at (0,0) [] {$ $};
	\node[vertex,label=270:{ $y$}] (y) at (4,0) [] {$ $};
	\node[vertex,label=270:{ $z$}] (z) at (8,0) [] {$ $};
	\node[vertex,label=90:{ $a$}] (a) at (2,7) [] {$ $};
	\draw (2,4) node {$f$};
	\node (p1) at(-4,0) [] {$ $};
	\node (p2) at(8,0) [] {$ $};
	\node (q1) at(-8,1) [] {$ $};
	\node[label=0:{ $C$}] (q2) at(12,1) [] {$ $};
	\path[-,out=-30, in=180] (q1) edge (p1);
	\draw (v) -- (z);
	\path[-,out=0, in=-150] (p2) edge (q2);
	\draw (v) -- (a);
	\draw (z) -- (a);
	\end{tikzpicture}
\end{center}
\caption{Case~3}\label{fig:3face}
\end{figure}

Let $f_1$, $f_2$ and $f_3$ be the (possibly identical) opposite faces of $f$ with respect to the $C$-edges $vx,xy,yz$ of $f$ (see Figure~\ref{fig:3face}). Then $f_2$ is not a minor 2-face for the same reason as in Case~2c. We distinguish the following subcases.

\begin{description}
\item[Case 3a:] \textit{Neither $f_1$ nor $f_3$ is a minor $3$-face \emph{(see Figure~\ref{fig:3face22})}.}

Then $f_2$ is neither a minor 2-face nor a minor 3-face, and $f_1$ and $f_3$ are minor 2-faces, as otherwise by R1--R5 $f$ would not send a total weight of more than~\nicefrac{2}{3} to its opposite faces. Moreover, $b\neq d$ (in the notation of Figure~\ref{fig:3face22}), since $xy$ is not extendable. Then $\overline{C}$ is obtained from $C$ by replacing the path $(w,v,x,y,z,u)$ with the path $(w,b,x,v,a,z,y,d,u)$.

\begin{figure}[!htb]
\begin{minipage}[b]{.5\textwidth}
\begin{center}
	\begin{tikzpicture}[scale=0.25, -, 
	vertex/.style={circle,fill=black,draw,minimum size=5pt,inner sep=0pt}]
	\node[vertex,label=270:{ $w$}] (w) at (-8,0) [] {$ $};
	\node[vertex,label=270:{ $v$}] (v) at (-4,0) [] {$ $};
	\node[vertex,label=270:{ $x$}] (x) at (0,0) [] {$ $};
	\node[vertex,label=270:{ $y$}] (y) at (4,0) [] {$ $};
	\node[vertex,label=270:{ $z$}] (z) at (8,0) [] {$ $};
	\node[vertex,label=270:{ $u$}] (u) at (12,0) [] {$ $};
	\node[vertex,label=90:{ $a$}] (a) at (2,7) [] {$ $};
	\draw (2,4) node {$f$};
	\draw (-4,-4.5) node {$f_1$};
	\draw (8,-4.5) node {$f_3$};
	\node[vertex,label=-90:{ $b$}] (b) at (-4,-7) [] {$ $};
	\node[vertex,label=-90:{ $d$}] (d) at (8,-7) [] {$ $};
	\node (p1) at(-8,0) [] {$ $};
	\node (p2) at(12,0) [] {$ $};
	\node (q1) at(-12,1) [] {$ $};
	\node[label=0:{ $C$}] (q2) at(16,1) [] {$ $};
	\path[-,out=-30, in=180] (q1) edge (p1);
	\draw (w) -- (u);
	\path[-,out=0, in=-150] (p2) edge (q2);
	\draw (v) -- (a);
	\draw (z) -- (a);
	\draw (w) -- (b);
	\draw (x) -- (b);
	\draw (y) -- (d);
	\draw (u) -- (d);
	\end{tikzpicture}
\caption{Case~3a}\label{fig:3face22}
\end{center}
\end{minipage}\hfil %
\begin{minipage}[b]{.5\textwidth}
\begin{center}
	\begin{tikzpicture}[scale=0.25, -, 
	vertex/.style={circle,fill=black,draw,minimum size=5pt,inner sep=0pt}]
	\node[vertex,label=270:{ $v$}] (v) at (-4,0) [] {$ $};
	\node[vertex,label=270:{ $x$}] (x) at (0,0) [] {$ $};
	\node[vertex,label=270:{ $y$}] (y) at (4,0) [] {$ $};
	\node[vertex,label=270:{ $z$}] (z) at (8,0) [] {$ $};
	\node[vertex,label=270:{ $u$}] (u) at (12,0) [] {$ $};
	\node[vertex,label=90:{ $a$}] (a) at (2,7) [] {$ $};
	\draw (2,4) node {$f$};
	\draw (6,-4.5) node {$f_3$};
	\node[vertex,label=-90:{ $d$}] (d) at (6,-7) [] {$ $};
	\node (p1) at(-4,0) [] {$ $};
	\node (p2) at(12,0) [] {$ $};
	\node (q1) at(-8,1) [] {$ $};
	\node[label=0:{ $C$}] (q2) at(16,1) [] {$ $};
	\path[-,out=-30, in=180] (q1) edge (p1);
	\draw (v) -- (u);
	\path[-,out=0, in=-150] (p2) edge (q2);
	\draw (v) -- (a);
	\draw (z) -- (a);
	\draw (x) -- (d);
	\draw (u) -- (d);
	\end{tikzpicture}
\caption{Case~3b}\label{fig:3face23}
\end{center}
\end{minipage}
\end{figure}

\item[Case 3b:] \textit{$f_1$ or $f_3$ is a minor $3$-face \emph{(see Figure~\ref{fig:3face23})}.}

The face $f_2$ is not a minor 3-face with middle $C$-edge $xy$, as otherwise $\{v,z\}$ would be a 2-separator of $G$. Hence, $f_1 \neq f_3$. Since $f$ sends a total weight of more than \nicefrac{2}{3} to its opposite faces, at least one of $f_1$ and $f_3$ is a minor 3-face that has its middle $C$-edge in $\{vx,yz\}$ by~R3, say without loss of generality that the middle $C$-edge of $f_3$ is $yz$. Then $\overline{C}$ is obtained from $C$ by replacing the path $(v,x,y,z,u)$ with the path $(v,a,z,y,x,d,u)$.
\end{description}

\item[Case 4:] \textit{$f$ is a minor $4$-face \emph{(see Figure~\ref{fig:4face})}.}

Then $f$ is initially charged with weight 4. If $f$ looses a total net weight of at most \nicefrac{2}{3}, then $w(f) \geq \nicefrac{10}{3}$, so assume that weight at least \nicefrac{3}{3} is sent to opposite faces. We have to show that this is impossible by considering Rules R2--R5.

\begin{figure}[!htb]
\begin{center}
	\begin{tikzpicture}[scale=0.25, -, 
	vertex/.style={circle,fill=black,draw,minimum size=5pt,inner sep=0pt}]
	\node[vertex,label=270:{ $v$}] (v) at (-8,0) [] {$ $};
	\node[vertex,label=270:{ $w$}] (w) at (-4,0) [] {$ $};
	\node[vertex,label=270:{ $x$}] (x) at (0,0) [] {$ $};
	\node[vertex,label=270:{ $y$}] (y) at (4,0) [] {$ $};
	\node[vertex,label=270:{ $z$}] (z) at (8,0) [] {$ $};
	\node[vertex,label=90:{ $a$}] (a) at (0,7) [] {$ $};
	\draw (0,4) node {$f$};
	\node (p1) at(-8,0) [] {$ $};
	\node (p2) at(8,0) [] {$ $};
	\node (q1) at(-12,1) [] {$ $};
	\node[label=0:{ $C$}] (q2) at(12,1) [] {$ $};
	\path[-,out=-30, in=180] (q1) edge (p1);
	\draw (v) -- (z);
	\path[-,out=0, in=-150] (p2) edge (q2);
	\draw (v) -- (a);
	\draw (z) -- (a);
	\end{tikzpicture}
\end{center}
\caption{Case~4}\label{fig:4face}
\end{figure}

Assume first that $f$ has an opposite minor 2-face $f'$. We distinguish the following subcases.

\begin{description}
\item[Case 4a:] \textit{$f'$ has $C$-edges $wx$ and $xy$ \emph{(see Figure~\ref{fig:4face22a})}.}

Then $vx$ or $xz$ is an edge of $G$ and $C$ can be extended by detouring $C$ through one of these edges and $d$, which contradicts the choice of $C$.

\begin{figure}[!htb]
\begin{minipage}[b]{.5\textwidth}
\begin{center}
	\begin{tikzpicture}[scale=0.25, -, 
	vertex/.style={circle,fill=black,draw,minimum size=5pt,inner sep=0pt}]
	\node[vertex,label=270:{ }] (f) at (-12,0) [] {$ $};
	\node[vertex,label=270:{ $v$}] (v) at (-8,0) [] {$ $};
	\node[vertex,label=270:{ $w$}] (w) at (-4,0) [] {$ $};
	\node[vertex,label=270:{ $x$}] (x) at (0,0) [] {$ $};
	\node[vertex,label=270:{ $y$}] (y) at (4,0) [] {$ $};
	\node[vertex,label=270:{ $z$}] (z) at (8,0) [] {$ $};
	\node[vertex,label=90:{ $a$}] (a) at (0,7) [] {$ $};
	\draw (0,4) node {$f$};
	\draw (0,-4) node {$f'$};
	\node[vertex,label=-90:{ $d$}] (d) at (0,-7) [] {$ $};
	\node (p1) at(-12,0) [] {$ $};
	\node (p2) at(8,0) [] {$ $};
	\node (q1) at(-16,1) [] {$ $};
	\node[label=0:{ $C$}] (q2) at(12,1) [] {$ $};
	\path[-,out=-30, in=180] (q1) edge (p1);
	\draw (f) -- (z);
	\path[-,out=0, in=-150] (p2) edge (q2);
	\path[dashed,out=45, in=135] (v) edge (x);
	\path[dashed,out=45, in=135] (x) edge (z);
	\draw (v) -- (a);
	\draw (z) -- (a);
	\draw (w) -- (d);
	\draw (y) -- (d);
	\end{tikzpicture}
\caption{Case~4a}\label{fig:4face22a}
\end{center}
\end{minipage}\hfil %
\begin{minipage}[b]{.5\textwidth}
\begin{center}
	\begin{tikzpicture}[scale=0.25, -, 
	vertex/.style={circle,fill=black,draw,minimum size=5pt,inner sep=0pt}]
	\node[vertex,label=270:{ $t$}] (t) at (-12,0) [] {$ $};
	\node[vertex,label=270:{ $v$}] (w) at (-8,0) [] {$ $};
	\node[vertex,label=270:{ $w$}] (v) at (-4,0) [] {$ $};
	\node[vertex,label=270:{ $x$}] (x) at (0,0) [] {$ $};
	\node[vertex,label=270:{ $y$}] (y) at (4,0) [] {$ $};
	\node[vertex,label=270:{ $z$}] (z) at (8,0) [] {$ $};
	\node[vertex,label=270:{ $u$}] (u) at (12,0) [] {$ $};
	\node[vertex,label=90:{ $a$}] (a) at (0,7) [] {$ $};
	\draw (0,4) node {$f$};
	\draw (8,-4) node {$f'$};
	\draw (-6,-4) node {$f_1$};
	\node[vertex,label=-90:{ $b$}] (b) at (-6,-7) [] {$ $};
	\node[vertex,label=-90:{ $d$}] (d) at (8,-7) [] {$ $};
	\node (p1) at(-12,0) [] {$ $};
	\node (p2) at(12,0) [] {$ $};
	\node (q1) at(-16,1) [] {$ $};
	\node[label=0:{ $C$}] (q2) at(16,1) [] {$ $};
	\path[-,out=-30, in=180] (q1) edge (p1);
	\draw (t) -- (u);
	\path[-,out=0, in=-150] (p2) edge (q2);
	\draw (w) -- (a);
	\draw (z) -- (a);
	\draw (t) -- (b);
	\draw (x) -- (b);
	\draw (y) -- (d);
	\draw (u) -- (d);
	\end{tikzpicture}
\end{center}
\caption{Case~4b}\label{fig:4face22aa}
\end{minipage}\hfil %
\end{figure}

\item[Case 4b:] \textit{Every opposite minor 2-face of $f$ has exactly one $C$-edge of $f$ \emph{(see Figure~\ref{fig:4face22aa})}.}

In particular, $m_{f,f'} = 1$. Without loss of generality, let $f'$ have the $C$-edge $yz$. Then $f$ sends weight \nicefrac{2}{3} to $f'$ by~R2, and R1 does not decrease the weight of $f$. 
Moreover, if $f$ sends weight to another face with the Rules~R4 or~R5, then $xy$ is a $C$-edge of a major face (since $C$ does not contain any extendable edge) and $f$ gains weight 1 from this major face, so that $w(f) \geq 4-\nicefrac{2}{3}+1-\nicefrac{2}{3} = \nicefrac{11}{3}$, which contradicts $w(f) < \nicefrac{10}{3}$. Therefore, $f$ has by~R2 and~R3 an opposite minor 2- or 3-face $f_1 \neq f'$. If $f_1$ is a minor 2-face, $m_{f,f_1} = 1$, so that $f_1$ has the $C$-edge $vw$. Then neither $wx$ nor $xy$ is a $C$-edge of a minor face opposite to $f$, as such a minor face would be a 2-face with $C$-edges $wx$ and $xy$ (see Case~4a). Thus, $f$ gains weight 2 from the major face(s) with $C$-edges $wx$ and $xy$, which contradicts $w(f) < \nicefrac{10}{3}$.

Hence, $f_1$ is a minor 3-face. Since $w(f) < \nicefrac{10}{3}$, the middle $C$-edge of $f_1$ is either $vw$ or $wx$. If it is $vw$, $\overline{C}$ can be obtained from $C$ by replacing the path $(t,v,w,x,y,z,u)$ with $(t,b,x,w,v,a,z,y,d,u)$ (see Figure~\ref{fig:4face22aa}), as we have $b \neq d$, since otherwise $C$ would contain the extendable edge $xy$. Hence, let the middle $C$-edge of $f_1$ be $wx$. Then $wz \notin E(G)$, as otherwise $C$ could be extended by replacing the path $(v,w,x,y,z)$ with $(v,b,y,x,w,z)$. Since $\{v,y\}$ is not a 2-separator of the 3-connected graph $G$, this implies $xz \in E(G)$. Then $\overline{C}$ can be obtained from $C$ by replacing the path $(x,y,z,u)$ with $(x,z,y,d,u)$, which contradicts the choice of $C$.
\end{description}

\item From Cases 4a+b, we conclude that $f'$ has either the $C$-edges $vw$ and $wx$ or $C$-edges $xy$ and $yz$, say without loss of generality the latter.

\begin{description}
\item[Case 4c:] \textit{$f'$ has $C$-edges $xy$ and $yz$, and $f$ has an opposite major face \emph{(see Figure~\ref{fig:4face2f2a})}.}

Then $wy \notin E(G)$, as otherwise $C$ can be extended by detouring through $f'$. Hence, $vy \in E(G)$, as otherwise $\deg_G(y)=2$. Since $f$ has an opposite major face and $wx$ is not an extendable edge of $C$, $wx$ is a $C$-edge of such an opposite major face $f''$. Then $f$ gains weight 1 from $f''$ by~R1 and sends by R2 weight \nicefrac{2}{3} to a minor opposite 2-face with $C$-edge $vw$ in order to satisfy the assumption $w(f) < \nicefrac{10}{3}$ (see Figure~\ref{fig:4face2f2a} and note that~R4 and~R5 do not apply here). But this is impossible, as then $C$ can be extended by replacing the path $(t,v,w,x,y,z)$ with $(t,b,w,v,y,x,d,z)$, since $b \neq d$.

\begin{figure}[!htb]
\begin{minipage}[b]{.5\textwidth}
\begin{center}
	\begin{tikzpicture}[scale=0.25, -, 
	vertex/.style={circle,fill=black,draw,minimum size=5pt,inner sep=0pt}]
	\node[vertex,label=270:{ $t$}] (f) at (-12,0) [] {$ $};
	\node[vertex,label=270:{ $v$}] (v) at (-8,0) [] {$ $};
	\node[vertex,label=270:{ $w$}] (w) at (-4,0) [] {$ $};
	\node[vertex,label=270:{ $x$}] (x) at (0,0) [] {$ $};
	\node[vertex,label=270:{ $y$}] (y) at (4,0) [] {$ $};
	\node[vertex,label=270:{ $z$}] (z) at (8,0) [] {$ $};
	\node[vertex,label=90:{ $a$}] (a) at (0,7) [] {$ $};
	\draw (0,4) node {$f$};
	\draw (4,-4) node {$f'$};
	\node[vertex,label=-90:{ $b$}] (b) at (-8,-7) [] {$ $};
	\node[vertex,label=-90:{ $d$}] (d) at (4,-7) [] {$ $};
	\node (p1) at(-12,0) [] {$ $};
	\node (p2) at(8,0) [] {$ $};
	\node (q1) at(-16,1) [] {$ $};
	\node[label=0:{ $C$}] (q2) at(12,1) [] {$ $};
	\path[-,out=-30, in=180] (q1) edge (p1);
	\draw (f) -- (z);
	\path[-,out=0, in=-150] (p2) edge (q2);
	\path[dashed,out=45, in=135] (v) edge (y);
	\draw (v) -- (a);
	\draw (z) -- (a);
	\draw (f) -- (b);
	\draw (w) -- (b);
	\draw (x) -- (d);
	\draw (z) -- (d);
	\end{tikzpicture}
\caption{Case~4c}\label{fig:4face2f2a}
\end{center}
\end{minipage}
\begin{minipage}[b]{.5\textwidth}
\begin{center}
	\begin{tikzpicture}[scale=0.25, -, 
	vertex/.style={circle,fill=black,draw,minimum size=5pt,inner sep=0pt}]
	\node[vertex,label=270:{ $t$}] (v) at (-8,0) [] {$ $};
	\node[vertex,label=270:{ $v$}] (w) at (-4,0) [] {$ $};
	\node[vertex,label=270:{ $w$}] (x) at (0,0) [] {$ $};
	\node[vertex,label=270:{ $x$}] (y) at (4,0) [] {$ $};
	\node[vertex,label=270:{ $y$}] (z) at (8,0) [] {$ $};
	\node[vertex,label=270:{ $z$}] (t) at (12,0) [] {$ $};
	\node[vertex,label=90:{ $a$}] (a) at (4,7) [] {$ $};
	\draw (4,4) node {$f$};
	\draw (-2,-4) node {$f_1$};
	\draw (8,-4) node {$f'$};
	\node[vertex,label=-90:{ $b$}] (b) at (-2,-7) [] {$ $};
	\node[vertex,label=-90:{ $d$}] (d) at (8,-7) [] {$ $};
	\node (p1) at(-8,0) [] {$ $};
	\node (p2) at(12,0) [] {$ $};
	\node (q1) at(-12,1) [] {$ $};
	\node[label=0:{ $C$}] (q2) at(16,1) [] {$ $};
	\path[-,out=-30, in=180] (q1) edge (p1);
	\draw (v) -- (t);
	\path[-,out=0, in=-150] (p2) edge (q2);
	\path[dashed,out=45, in=135] (w) edge (z);
	\draw (w) -- (a);
	\draw (t) -- (a);
	\draw (v) -- (b);
	\draw (y) -- (b);
	\draw (y) -- (d);
	\draw (t) -- (d);
	\end{tikzpicture}
\caption{Case~4d}\label{fig:4face22}
\end{center}
\end{minipage}\hfil %
\end{figure}

\item[Case 4d:] \textit{$f'$ has $C$-edges $xy$ and $yz$, and $wx$ is a $C$-edge of a minor 2- or 3-face $f_1$ \emph{(see Figure~\ref{fig:4face22})}.}

As in Case~4c, $wy \notin E(G)$ and $vy \in E(G)$. Hence, $f_1$ is a minor 3-face, as otherwise $\deg_G(w)=2$. Then $\overline{C}$ is obtained from $C$ by replacing the path $(t,v,w,x,y,z)$ with $(t,b,x,w,v,y,z)$ (note that $b=d$ is possible).

\item[Case 4e:] \textit{$f'$ has $C$-edges $xy$ and $yz$, and $wx$ is a $C$-edge of a minor $j$-face $f_1$ with $j \geq 4$ \emph{(see Figure~\ref{fig:4face2f2})}.}

Then $f$ gains weight \nicefrac{2}{3} from $f_1$ by~R4 and sends weight \nicefrac{4}{3} to $f'$. Hence, we get the contradiction $w(f) = \nicefrac{10}{3}$, unless $f$ sends weight \nicefrac{2}{3} to $f_1$ by~R4 or \nicefrac{1}{3} to $f_1$ by~R5. In that case, $j = 4$ or $j = 5$ and there are only  minor 2-faces opposite to $f_1$. 
As argued in Case~4c, $wy \notin E(G)$ and $vy \in E(G)$. Moreover, $uw$ (and $su$ in case of $j=5$; see Figure~\ref{fig:4face2f2}) are not edges of $G$, as otherwise $C$ can be extended by detouring through $g$. Hence, $ux \in E(G)$, as otherwise $\deg_G(u)=2$, which is a contradiction. This implies $\deg_G(w)=2$, which is a contradiction.

\begin{figure}[!htb]
\begin{center}
	\begin{tikzpicture}[scale=0.25, -, 
	vertex/.style={circle,fill=black,draw,minimum size=5pt,inner sep=0pt}]
	\node[vertex,label=270:{ $r$}] (r) at (-20,0) [] {$ $};
	\node[vertex,label=270:{ $s$}] (s) at (-16,0) [] {$ $};
	\node[vertex,label=270:{ $t$}] (t) at (-12,0) [] {$ $};
	\node[vertex,label=270:{ $u$}] (u) at (-8,0) [] {$ $};
	\node[vertex,label=270:{ $v$}] (v) at (-4,0) [] {$ $};
	\node[vertex,label=270:{ $w$}] (w) at (0,0) [] {$ $};
	\node[vertex,label=270:{ $x$}] (x) at (4,0) [] {$ $};
	\node[vertex,label=270:{ $y$}] (y) at (8,0) [] {$ $};
	\node[vertex,label=270:{ $z$}] (z) at (12,0) [] {$ $};
	\draw (4,4) node {$f$};
	\draw (-6,-4) node {$f_1$};
	\draw (-16,4) node {$f_2$};
	\draw (-8,4) node {$f_3$};
	\draw (8,-4) node {$f'$};
	\node[vertex,label=90:{ $a$}] (a) at (4,7) [] {$ $};
	\node[vertex,label=90:{ $g$}] (g) at (-8,7) [] {$ $};
	\node[vertex,label=90:{ $h$}] (h) at (-16,7) [] {$ $};
	\node[vertex,label=-90:{ $b$}] (b) at (-6,-7) [] {$ $};
	\node[vertex,label=-90:{ $d$}] (d) at (8,-7) [] {$ $};
	\node (p1) at(-20,0) [] {$ $};
	\node (p2) at(12,0) [] {$ $};
	\node (q1) at(-24,1) [] {$ $};
	\node[label=0:{ $C$}] (q2) at(16,1) [] {$ $};
	\path[-,out=-30, in=180] (q1) edge (p1);
	\draw (r) -- (z);
	\path[-,out=0, in=-150] (p2) edge (q2);
	\path[dashed,out=45, in=135] (v) edge (y);
	\path[dashed,out=-35, in=-145] (u) edge (x);
	\draw (v) -- (a);
	\draw (z) -- (a);
	\draw (s) -- (b);
	\draw (x) -- (b);
	\draw (x) -- (d);
	\draw (z) -- (d);
	\draw (r) -- (h);
	\draw (t) -- (h);
	\draw (t) -- (g);
	\draw (v) -- (g);
	\end{tikzpicture}
\caption{Case~4e}\label{fig:4face2f2}
\end{center}
\end{figure}
\end{description}

From Cases~4a--e, we conclude that $f$ has no opposite minor 2-face. Then $w(f)<\nicefrac{10}{3}$ and~R1--R5 imply that $f$ has an opposite minor 3-face that has a $C$-edge of $f$ as middle $C$-edge (due to R3), or an opposite minor 4-face $f'$ with $m_{f,f'}=2$ that has an opposite minor 2- or 3-face $f_2$ with $m_{f',f_2}=2$ (due to R4); note that we still contradict $w(f) <\nicefrac{10}{3}$ when $f$ has two opposite minor 5-faces, to each of which $f$ sends weight \nicefrac{1}{3} by~R5. We therefore distinguish these remaining subcases.

\begin{description}
\item[Case 4f:] \textit{$f$ has an opposite minor 3-face $f'$ with middle $C$-edge $wx$ or $xy$ \emph{(see Figure~\ref{fig:4face323})}.}

Without loss of generality, let $xy$ be the middle $C$-edge of $f'$. Then $vy \notin E(G)$, as otherwise $C$ can be extended by replacing the path $(v,w,x,y,z)$ with $(v,y,x,w,d,z)$. This implies $wy \in E(G)$, as otherwise $\deg_G(y) = 2$. Since $\{w,z\}$ is no 2-separator of $G$, $vx \in E(G)$. Then $C$ can be extended by replacing the path $(v,w,x,y,z)$ with $(v,x,y,w,d,z)$.

\begin{figure}[!htb]
\begin{minipage}[b]{.5\textwidth}
\begin{center}
	\begin{tikzpicture}[scale=0.25, -, 
	vertex/.style={circle,fill=black,draw,minimum size=5pt,inner sep=0pt}]
	\node[vertex,label=270:{ $v$}] (v) at (-4,0) [] {$ $};
	\node[vertex,label=270:{ $w$}] (x) at (0,0) [] {$ $};
	\node[vertex,label=270:{ $x$}] (y) at (4,0) [] {$ $};
	\node[vertex,label=270:{ $y$}] (z) at (8,0) [] {$ $};
	\node[vertex,label=270:{ $z$}] (u) at (12,0) [] {$ $};
	\node[vertex,label=90:{ $a$}] (a) at (4,7) [] {$ $};
	\draw (4,4) node {$f$};
	\draw (6,-4.5) node {$f'$};
	\node[vertex,label=-90:{ $d$}] (d) at (6,-7) [] {$ $};
	\node (p1) at(-4,0) [] {$ $};
	\node (p2) at(12,0) [] {$ $};
	\node (q1) at(-8,1) [] {$ $};
	\node[label=0:{ $C$}] (q2) at(16,1) [] {$ $};
	\path[-,out=-30, in=180] (q1) edge (p1);
	\path[dashed,out=45, in=135] (v) edge (y);
	\draw (v) -- (u);
	\path[-,out=0, in=-150] (p2) edge (q2);
	\path[dashed,out=-45, in=-135] (x) edge (z);
	\draw (v) -- (a);
	\draw (u) -- (a);
	\draw (x) -- (d);
	\draw (u) -- (d);
	\end{tikzpicture}
\caption{Case~4f}\label{fig:4face323}
\end{center}
\end{minipage}\hfil %
\begin{minipage}[b]{.5\textwidth}
\begin{center}
	\begin{tikzpicture}[scale=0.25, -, 
	vertex/.style={circle,fill=black,draw,minimum size=5pt,inner sep=0pt}]
	\node[vertex,label=270:{ $t$}] (t) at (-12,0) [] {$ $};
	\node[vertex,label=270:{ $v$}] (v) at (-8,0) [] {$ $};
	\node[vertex,label=270:{ $w$}] (w) at (-4,0) [] {$ $};
	\node[vertex,label=270:{ $x$}] (x) at (0,0) [] {$ $};
	\node[vertex,label=270:{ $y$}] (y) at (4,0) [] {$ $};
	\node[vertex,label=270:{ $z$}] (z) at (8,0) [] {$ $};
	\node[vertex,label=270:{ $u$}] (u) at (12,0) [] {$ $};
	\node[vertex,label=90:{ $a$}] (a) at (0,7) [] {$ $};
	\draw (0,4) node {$f$};
	\draw (6,-4) node {$f'$};
	\draw (-6,-4) node {$f_1$};
	\node[vertex,label=-90:{ $b$}] (b) at (-6,-7) [] {$ $};
	\node[vertex,label=-90:{ $d$}] (d) at (6,-7) [] {$ $};
	\node (p1) at(-12,0) [] {$ $};
	\node (p2) at(12,0) [] {$ $};
	\node (q1) at(-16,1) [] {$ $};
	\node[label=0:{ $C$}] (q2) at(16,1) [] {$ $};
	\path[-,out=-30, in=180] (q1) edge (p1);
	\draw (t) -- (u);
	\path[-,out=0, in=-150] (p2) edge (q2);
	\path[dashed,out=-45, in=-135] (t) edge (w);
	\path[dashed,out=-45, in=-135] (y) edge (u);
	\draw (v) -- (a);
	\draw (z) -- (a);
	\draw (t) -- (b);
	\draw (x) -- (b);
	\draw (x) -- (d);
	\draw (u) -- (d);
	\end{tikzpicture}
\caption{Case~4g}\label{fig:4face422}
\end{center}
\end{minipage}
\end{figure}

\item[Case 4g:] \textit{$f$ has an opposite minor 3-face $f'$ with middle $C$-edge $vw$ or $yz$, but no opposite 4-face \emph{(see Figure~\ref{fig:4face422})}.}

Without loss of generality, let $yz$ be the middle $C$-edge of $f'$. Let $f_1$ be the face opposite to $f$ that has $C$-edge $wx$. Then $f_1$ is not major, as otherwise $w(f) = 4-1+1 > \nicefrac{10}{3}$, since $f$ has no opposite minor 2-faces. For the same reason, $f_1$ is a minor $j$-face satisfying $j \geq 3$. If $j \geq 5$, $f_1$ sends weight $\nicefrac{2}{3}$ to $f$ due to~R4, which contradicts $w(f) < \nicefrac{10}{3}$, as $f$ sends weight at most \nicefrac{1}{3} to $f_1$ due to~R5 (exactly \nicefrac{1}{3} only if $j=5$ and $f_1$ has two opposite 2-faces).

Since $j \neq 4$ by assumption, $f_1$ is a minor 3-face (see Figure~\ref{fig:4face422}). Then $wy \notin E(G)$, as otherwise $\overline{C}$ is obtained from $C$ by replacing the path $(v,w,x,y,z,u)$ with $(v,a,z,y,w,x,d,u)$, and $wz \notin E(G)$, as otherwise $\overline{C}$ is obtained from $C$ by replacing the path $(w,x,y,z,u)$ with $(w,z,y,x,d,u)$. Hence, $tw \in E(G)$, as otherwise $\deg_G(w) = 2$.
Then $\overline{C}$ is obtained from $C$ by replacing the path $(t,v,w,x,y,z,u)$ with $(t,w,v,a,z,y,x,d,u)$, which contradicts the choice of $C$.

\item[Case 4h:] \textit{$f$ has an opposite minor 3-face $f'$ with middle $C$-edge $vw$ or $yz$ and an opposite 4-face $f_1$ \emph{(see Figure~\ref{fig:4face327})}.}

Without loss of generality, let $yz$ be the middle $C$-edge of $f'$. Then $m_{f,f_1}=2$, as otherwise $wx$ is a $C$-edge of a major face, which would imply $w(f) =4-1+1 > \nicefrac{10}{3}$. Hence, $f_1$ sends weight \nicefrac{2}{3} to $f$ by~R4, which implies that $f$ must send weight \nicefrac{2}{3} to $f_1$ by~R4, as otherwise $w(f) \geq \nicefrac{10}{3}$. Hence, $f_1$ has an opposite minor 2- or 3-face $f_2$ that satisfies $m_{f_1,f_2} = 2$ (see Figure~\ref{fig:4face327}). Then $wy \notin E(G)$, as otherwise $C$ can be extended by replacing the path $(v,w,x,y,z,q)$ with $(v,a,z,y,w,x,d,q)$, and $wz \notin E(G)$, as otherwise $C$ can be extended by replacing the path $(w,x,y,z,q)$ with $(w,z,y,x,d,q)$. If $f_2$ is a 3-face, this implies by symmetry $tw \notin E(G)$ and $uw \notin E(G)$, which contradicts $\deg_G(w) \geq 3$. Hence, $f_2$ is a 2-face. Then $uw \notin E(G)$, as otherwise $C$ can be extended by replacing the path $(t,u,v,w)$ with $(t,g,v,u,w)$, which implies $tw \in E(G)$, as otherwise $\deg_G(w)=2$. This contradicts $\deg_G(u) \geq 3$.

\begin{figure}[!htb]
\begin{center}
	\begin{tikzpicture}[scale=0.25, -, 
	vertex/.style={circle,fill=black,draw,minimum size=5pt,inner sep=0pt}]
	\node[vertex,label=270:{ $t$}] (r) at (-20,0) [] {$ $};
	\node[vertex,label=270:{ $u$}] (s) at (-16,0) [] {$ $};
	\node[vertex,label=270:{ $v$}] (t) at (-12,0) [] {$ $};
	\node[vertex,label=270:{ $w$}] (u) at (-8,0) [] {$ $};
	\node[vertex,label=270:{ $x$}] (v) at (-4,0) [] {$ $};
	\node[vertex,label=270:{ $y$}] (w) at (0,0) [] {$ $};
	\node[vertex,label=270:{ $z$}] (x) at (4,0) [] {$ $};
	\node[vertex,label=270:{ $q$}] (y) at (8,0) [] {$ $};
	\node[vertex,label=-90:{ $b$}] (b) at (-12,-7) [] {$ $};
	\node[vertex,label=-90:{ $d$}] (d) at (2,-7) [] {$ $};
	\draw (2,-4) node {$f'$};
	\draw (-12,-4) node {$f_1$};
	\draw (-4,4) node {$f$};
	\draw (-16,4) node {$f_2$};
	\node[vertex,label=90:{ $a$}] (a) at (-4,7) [] {$ $};
	\node[vertex,label=90:{ $g$}] (g) at (-16,7) [] {$ $};
	\node (p1) at(-20,0) [] {$ $};
	\node (p2) at(8,0) [] {$ $};
	\node (q1) at(-24,1) [] {$ $};
	\node[label=0:{ $C$}] (q2) at(12,1) [] {$ $};
	\path[-,out=-30, in=180] (q1) edge (p1);
	\draw (r) -- (y);
	\path[-,out=0, in=-150] (p2) edge (q2);
	\path[dashed,out=-45, in=-135] (r) edge (u);
	\draw (r) -- (b);
	\draw (v) -- (b);
	\draw (v) -- (d);
	\draw (y) -- (d);
	\draw (t) -- (a);
	\draw (x) -- (a);
	\draw (r) -- (g);
	\draw (t) -- (g);
	\end{tikzpicture}
\caption{Case~4h}\label{fig:4face327}
\end{center}
\end{figure}

\item[Case 4i:] \textit{$f$ has no opposite minor 3-face whose middle $C$-edge is a $C$-edge of $f$ \emph{(see Figure~\ref{fig:4face32})}.}

Then, as argued before, $f$ has an opposite minor 4-face $f'$ with $m_{f,f'}=2$ and $C$-edges $xy$ and $yz$, that has an opposite minor 2- or 3-face $f_2$ with $m_{f',f_2}=2$. According to~R4, $f$ sends weight \nicefrac{2}{3} to $f'$. Let $f''$ be the face opposite to $f$ that has $C$-edge $wx$. Then $f''$ must be either a second opposite minor 4-face with $m_{f,f''}=2$ that has an opposite minor 2- or 3-face $f_1$ with $m_{f'',f_1}=2$ (due to~R4), or a opposite minor 5-face with $m_{f,f''}=2$ that has two opposite minor 2-faces (due to~R5), as otherwise $w(f) \geq 4-\nicefrac{2}{3} = \nicefrac{10}{3}$, since $f$ sends no weight to any 2- or 3-face by~R2 or~R3.
Note that $g=a=h$ and $b=d$ are possible.

\begin{figure}[!htb]
\begin{center}
	\begin{tikzpicture}[scale=0.25, -, 
	vertex/.style={circle,fill=black,draw,minimum size=5pt,inner sep=0pt}]
	\node[vertex,label=270:{ $s$}] (q) at (-24,0) [] {$ $};
	\node[vertex,label=270:{ $t$}] (r) at (-20,0) [] {$ $};
	\node[vertex,label=270:{ $u$}] (s) at (-16,0) [] {$ $};
	\node[vertex,label=270:{ $v$}] (t) at (-12,0) [] {$ $};
	\node[vertex,label=270:{ $w$}] (u) at (-8,0) [] {$ $};
	\node[vertex,label=270:{ $x$}] (v) at (-4,0) [] {$ $};
	\node[vertex,label=270:{ $y$}] (w) at (0,0) [] {$ $};
	\node[vertex,label=270:{ $z$}] (x) at (4,0) [] {$ $};
	\node[vertex,label=270:{ $q$}] (y) at (8,0) [] {$ $};
	\node[vertex,label=270:{ $r$}] (z) at (12,0) [] {$ $};
	\node[vertex,label=-90:{ $b$}] (b) at (-12,-7) [] {$ $};
	\node[vertex,label=-90:{ $d$}] (d) at (4,-7) [] {$ $};
	\draw (4,-4) node {$f'$};
	\draw (-12,-4) node {$f''$};
	\draw (-4,4) node {$f$};
	\draw (8,4) node {$f_2$};
	\draw (-18,4) node {$f_1$};
	\node[vertex,label=90:{ $a$}] (a) at (-4,7) [] {$ $};
	\node[vertex,label=90:{ $g$}] (g) at (-18,7) [] {$ $};
	\node[vertex,label=90:{ $h$}] (h) at (8,7) [] {$ $};
	\node (p1) at(-24,0) [] {$ $};
	\node (p2) at(12,0) [] {$ $};
	\node (q1) at(-28,1) [] {$ $};
	\node[label=0:{ $C$}] (q2) at(16,1) [] {$ $};
	\path[-,out=-30, in=180] (q1) edge (p1);
	\draw (q) -- (z);
	\path[-,out=0, in=-150] (p2) edge (q2);
	\path[dashed,out=-45, in=-135] (v) edge (y);
	\path[dashed,out=45, in=135] (t) edge (w);
	\draw (r) -- (b);
	\draw (v) -- (b);
	\draw (v) -- (d);
	\draw (z) -- (d);
	\draw (t) -- (a);
	\draw (x) -- (a);
	\draw (x) -- (h);
	\draw (z) -- (h);
	\draw (q) -- (g);
	\draw (t) -- (g);
	\end{tikzpicture}
\caption{Case~4i}\label{fig:4face32}
\end{center}
\end{figure}

We claim that in all cases $vy$ is an edge of $G$. Consider the case that $f_2$ is a 2-face (see Figure~\ref{fig:4face32}). Then $yq \notin E(G)$, as otherwise $C$ can be extended by replacing the path $(y,z,q,r)$ with $(y,q,z,h,r)$, and thus $xq \in E(G)$, as otherwise $\deg_G(q)=2$. This implies that $vy$ or $wy$ is in $G$, as otherwise $\deg_G(y)=2$. Since $wy \notin E(G)$, as otherwise $C$ can be extended by replacing the path $(w,x,y,z,q,r)$ with $(w,y,x,q,z,h,r)$, we have $vy \in E(G)$, as claimed. Now consider the remaining case that $f_2$ is a 3-face. By symmetry, we will assume instead that $f_1$ is a 3-face and prove that $wz \in E(G)$ (such that the notation of Figure~\ref{fig:4face32} can be used); this implies $vy \in E(G)$ for the case that $f_2$ is a 3-face. Then $wy \notin E(G)$, as otherwise $C$ can be extended by replacing the path $(s,t,u,v,w,x,y)$ with $(s,g,v,u,t,b,x,w,y)$, and $uw \notin E(G)$, as otherwise $C$ can be extended by replacing the path $(s,t,u,v,w,x)$ with $(s,g,v,w,u,t,b,x)$. In addition, $tw \notin E(G)$, as otherwise $C$ can be extended by replacing the path $(s,t,u,v,w)$ with $(s,g,v,u,t,w)$. Then $wz \in E(G)$, as claimed, since otherwise $\deg_G(w)=2$, which is a contradiction.

Hence, we proved that in all cases $vy \in E(G)$. If $f''$ is a 5-face, then $ux \in E(G)$ by the last argument of Case~4e, which contradicts $\deg_G(w) \geq 3$. Hence, $f''$ is a 4-face, and no matter whether $f_1$ is a 2- or 3-face, $wz$ is an edge of $G$ by a symmetric argument to the one of the last paragraph. This contradicts that $G$ is plane, because $vy \in E(G)$.
\end{description}

\item[Case 5:] \textit{$f$ is a minor $5$-face \emph{(see Figure~\ref{fig:5face})}.}

Then $f$ is initially charged with weight 5. If $f$ looses a total net weight of at most \nicefrac{5}{3}, then $w(f) \geq \nicefrac{10}{3}$, so assume otherwise. We distinguish the following subcases.

\begin{figure}[!htb]
\begin{center}
	\begin{tikzpicture}[scale=0.25, -, 
	vertex/.style={circle,fill=black,draw,minimum size=5pt,inner sep=0pt}]
	\node[vertex,label=270:{ $u$}] (u) at (-12,0) [] {$ $};
	\node[vertex,label=270:{ $v$}] (v) at (-8,0) [] {$ $};
	\node[vertex,label=270:{ $w$}] (w) at (-4,0) [] {$ $};
	\node[vertex,label=270:{ $x$}] (x) at (0,0) [] {$ $};
	\node[vertex,label=270:{ $y$}] (y) at (4,0) [] {$ $};
	\node[vertex,label=270:{ $z$}] (z) at (8,0) [] {$ $};
	\node[vertex,label=90:{ $a$}] (a) at (-2,7) [] {$ $};
	\draw (-2,4) node {$f$};
	\node (p1) at(-12,0) [] {$ $};
	\node (p2) at(8,0) [] {$ $};
	\node (q1) at(-16,1) [] {$ $};
	\node[label=0:{ $C$}] (q2) at(12,1) [] {$ $};
	\path[-,out=-30, in=180] (q1) edge (p1);
	\draw (u) -- (z);
	\path[-,out=0, in=-150] (p2) edge (q2);
	\draw (u) -- (a);
	\draw (z) -- (a);
	\end{tikzpicture}
\end{center}
\caption{Case~5}\label{fig:5face}
\end{figure}

\begin{description}
\item[Case 5a:] \textit{$f$ sends weight to an opposite minor 5-face $f'$ \emph{(see Figure~\ref{fig:5faceA})}.}

Without loss of generality, let $xy$ and $yz$ be $C$-edges of $f'$ by~R5. Then $f$ sends weight \nicefrac{1}{3} to $f'$, and $f'$ has two opposite minor 2-faces $f_1$ and $f_2$. Since $w(f) < \nicefrac{10}{3}$, $f$ does neither send weight to a second 5-face nor to a 4-face nor to a 3-face (as there may be at most one of each kind and, if so, no 2-face that receives weight from $f$).
This implies that the edge $uv$ is a $C$-edge of a minor 2-face $f_3$ opposite to $f$, and that $vw$ and $wx$ are the $C$-edges of a second minor 2-face $f_4$ opposite to $f$ (see Figure~\ref{fig:5faceA}). Then $f'$ sends weight \nicefrac{1}{3} back to $f$ by~R5, but $w(f) = 5-3 \cdot \frac23 = 3 < \nicefrac{10}{3}$ is still satisfied.

\begin{figure}[!htb]
\begin{center}
	\begin{tikzpicture}[scale=0.25, -, 
	vertex/.style={circle,fill=black,draw,minimum size=5pt,inner sep=0pt}]
	\node[vertex,label=270:{ $t$}] (t) at (-16,0) [] {$ $};
	\node[vertex,label=270:{ $u$}] (u) at (-12,0) [] {$ $};
	\node[vertex,label=270:{ $v$}] (v) at (-8,0) [] {$ $};
	\node[vertex,label=270:{ $w$}] (w) at (-4,0) [] {$ $};
	\node[vertex,label=270:{ $x$}] (x) at (0,0) [] {$ $};
	\node[vertex,label=270:{ $y$}] (y) at (4,0) [] {$ $};
	\node[vertex,label=270:{ $z$}] (z) at (8,0) [] {$ $};
	\node[vertex,label=270:{ $p$}] (p) at (12,0) [] {$ $};
	\node[vertex,label=270:{ $q$}] (q) at (16,0) [] {$ $};
	\node[vertex,label=270:{ $r$}] (r) at (20,0) [] {$ $};
	\node[vertex,label=270:{ $s$}] (s) at (24,0) [] {$ $};
	\node[vertex,label=90:{ $a$}] (a) at (-2,7) [] {$ $};
	\node[vertex,label=-90:{ $i$}] (i) at (-12,-7) [] {$ $};
	\node[vertex,label=-90:{ $b$}] (b) at (-4,-7) [] {$ $};
	\node[vertex,label=-90:{ $d$}] (d) at (10,-7) [] {$ $};
	\node[vertex,label=90:{ $g$}] (g) at (12,7) [] {$ $};
	\node[vertex,label=90:{ $h$}] (h) at (20,7) [] {$ $};
	\draw (-2,4) node {$f$};
	\draw (10,-4) node {$f'$};
	\draw (12,4) node {$f_1$};
	\draw (20,4) node {$f_2$};
	\draw (-12,-4) node {$f_3$};
	\draw (-4,-4) node {$f_4$};
	\node (p1) at(-16,0) [] {$ $};
	\node (p2) at(24,0) [] {$ $};
	\node (q1) at(-20,1) [] {$ $};
	\node[label=0:{ $C$}] (q2) at(28,1) [] {$ $};
	\path[-,out=-30, in=180] (q1) edge (p1);
	\draw (t) -- (s);
	\path[-,out=0, in=-150] (p2) edge (q2);
	\draw (u) -- (a);
	\draw (z) -- (a);
	\path[dashed,out=35, in=145] (w) edge (z);
	\path[dashed,out=35, in=145] (w) edge (y);
	\path[dashed,out=-35, in=-145] (x) edge (p);
	\draw (x) -- (d);
	\draw (d) -- (r);
	\draw (v) -- (b);
	\draw (b) -- (x);
	\draw (z) -- (g);
	\draw (q) -- (g);
	\draw (q) -- (h);
	\draw (s) -- (h);
	\draw (i) -- (t);
	\draw (i) -- (v);
	\end{tikzpicture}
\end{center}
\caption{Case~5a}\label{fig:5faceA}
\end{figure}

We have $yp \notin E(G)$ and $pr \notin E(G)$, as otherwise $C$ can be extended by detouring through $g$. Since $\deg_G(p) \geq 3$, $xp \in E(G)$.
By symmetry, $wz \in E(G)$, which implies $yw \in E(G)$. Then $C$ can be extended by replacing the path $(v,w,x,y)$ with $(v,b,x,w,y)$.

\item[Case 5b:] \textit{$f$ sends weight to an opposite minor 4-face $f'$ \emph{(see Figure~\ref{fig:5faceB})}.}

Without loss of generality, let $xy$ and $yz$ be $C$-edges of $f'$ by~R4. Assume first that $f$ sends weight to an opposite minor 3-face $f_1$. Then $f$ sends total weight \nicefrac{5}{3} to $f'$ and $f_1$, and the middle $C$-edge of $f_1$ is either $uv$ or $vw$. Both cases contradict $w(f) <\nicefrac{10}{3}$, since no further weight is sent.
The same argument gives a contradiction if $f$ sends weight to a minor 4-face different from $f'$.

\begin{figure}[!htb]
\begin{center}
	\begin{tikzpicture}[scale=0.25, -, 
	vertex/.style={circle,fill=black,draw,minimum size=5pt,inner sep=0pt}]
	\node[vertex,label=270:{ $u$}] (u) at (-12,0) [] {$ $};
	\node[vertex,label=270:{ $v$}] (v) at (-8,0) [] {$ $};
	\node[vertex,label=270:{ $w$}] (w) at (-4,0) [] {$ $};
	\node[vertex,label=270:{ $x$}] (x) at (0,0) [] {$ $};
	\node[vertex,label=270:{ $y$}] (y) at (4,0) [] {$ $};
	\node[vertex,label=270:{ $z$}] (z) at (8,0) [] {$ $};
	\node[vertex,label=270:{ $q$}] (q) at (12,0) [] {$ $};
	\node[vertex,label=270:{ $r$}] (r) at (16,0) [] {$ $};
	\node[vertex,label=90:{ $a$}] (a) at (-2,7) [] {$ $};
	\node[vertex,label=-90:{ $b$}] (b) at (-4,-7) [] {$ $};
	\node[vertex,label=-90:{ $d$}] (d) at (8,-7) [] {$ $};
	\node[vertex,label=90:{ $g$}] (g) at (12,7) [] {$ $};
	\draw (-2,4) node {$f$};
	\draw (8,-4) node {$f'$};
	\draw (12,4) node {$f''$};
	\draw (-4,-4) node {$f_1$};
	\node (p1) at(-12,0) [] {$ $};
	\node (p2) at(16,0) [] {$ $};
	\node (q1) at(-16,1) [] {$ $};
	\node[label=0:{ $C$}] (q2) at(20,1) [] {$ $};
	\path[-,out=-30, in=180] (q1) edge (p1);
	\draw (u) -- (r);
	\path[-,out=0, in=-150] (p2) edge (q2);
	\draw (u) -- (a);
	\draw (z) -- (a);
	\path[-,out=0, in=-150] (p2) edge (q2);
	\path[dashed,out=35, in=145] (w) edge (z);
	\draw (x) -- (d);
	\draw (d) -- (r);
	\draw (v) -- (b);
	\draw (b) -- (x);
	\draw (z) -- (g);
	\draw (g) -- (r);
	\end{tikzpicture}
\end{center}
\caption{Case~5b}\label{fig:5faceB}
\end{figure}

Hence, $f$ sends a total weight of at least \nicefrac{4}{3} to minor 2-faces, as~R2 sends only multiples of weight \nicefrac{2}{3}. This implies that $f$ has an opposite minor 2-face $f_1$ with $m_{f,f_1}=2$. If $f_1$ has $C$-edges $uv$ and $vw$, then $wx$ is again a $C$-edge of major face, which sends weight 1 to $f$ and thus contradicts $w(f) < \nicefrac{10}{3}$. Hence, $f_1$ has $C$-edges $vw$ and $wx$ (see Figure~\ref{fig:5faceB}). Then $uw$ and $wy$ are not edges of $G$, as otherwise $C$ can be extended by detouring through $b$. Hence, $wz \in E(G)$, as otherwise $\deg_G(w)=2$.
Moreover, $yq \notin E(G)$ and $xq \in E(G)$ for the same reason as in Case~4i, which contradicts $\deg_G(y) \geq 3$.

\item[Case 5c:] \textit{$f$ sends weight to an opposite minor 3-face $f'$ with middle $C$-edge $wx$ \emph{(see Figure~\ref{fig:5faceC})}.}

In order to have $w(f) < \nicefrac{10}{3}$, by~R1--R3, $f$ sends weight \nicefrac{2}{3} to each of the minor 2-faces $f_1$ and $f_2$ having $C$-edges $uv$ and $yz$, respectively. Then $uw$ and $xz$ are not edges of $G$, as otherwise $C$ can be extended by detouring $C$ through $b$ or $g$, respectively. Since $\{v,y\}$ is not a 2-separator of $G$, this implies that either $wz \in E(G)$ or $ux \in E(G)$, say by symmetry the former. Then we can obtain $\overline{C}$ from $C$ by replacing the path $(v,w,x,y,z)$ with $(v,d,y,x,w,z)$.

\begin{figure}[!htb]
\begin{center}
	\begin{tikzpicture}[scale=0.25, -, 
	vertex/.style={circle,fill=black,draw,minimum size=5pt,inner sep=0pt}]
	\node[vertex,label=270:{ $r$}] (u) at (-12,0) [] {$ $};
	\node[vertex,label=270:{ $u$}] (v) at (-8,0) [] {$ $};
	\node[vertex,label=270:{ $v$}] (w) at (-4,0) [] {$ $};
	\node[vertex,label=270:{ $w$}] (x) at (0,0) [] {$ $};
	\node[vertex,label=270:{ $x$}] (y) at (4,0) [] {$ $};
	\node[vertex,label=270:{ $y$}] (z) at (8,0) [] {$ $};
	\node[vertex,label=270:{ $z$}] (q) at (12,0) [] {$ $};
	\node[vertex,label=270:{ $q$}] (r) at (16,0) [] {$ $};
	\node[vertex,label=90:{ $a$}] (a) at (2,7) [] {$ $};
	\node[vertex,label=-90:{ $b$}] (b) at (-8,-7) [] {$ $};
	\node[vertex,label=-90:{ $d$}] (d) at (2,-7) [] {$ $};
	\node[vertex,label=-90:{ $g$}] (g) at (12,-7) [] {$ $};
	\draw (2,4) node {$f$};
	\draw (2,-4) node {$f'$};
	\draw (-8,-4) node {$f_1$};
	\draw (12,-4) node {$f_2$};
	\node (p1) at(-12,0) [] {$ $};
	\node (p2) at(16,0) [] {$ $};
	\node (q1) at(-16,1) [] {$ $};
	\node[label=0:{ $C$}] (q2) at(20,1) [] {$ $};
	\path[-,out=-30, in=180] (q1) edge (p1);
	\draw (u) -- (r);
	\path[-,out=0, in=-150] (p2) edge (q2);
	\draw (v) -- (a);
	\draw (q) -- (a);
	\path[-,out=0, in=-150] (p2) edge (q2);
	\path[dashed,out=35, in=145] (x) edge (q);
	\draw (d) -- (z);
	\draw (d) -- (w);
	\draw (u) -- (b);
	\draw (b) -- (w);
	\draw (z) -- (g);
	\draw (g) -- (r);
	\end{tikzpicture}
\end{center}
\caption{Case~5c}\label{fig:5faceC}
\end{figure}

\item[Case 5d:] \textit{$f$ sends weight to an opposite minor 3-face $f'$ with middle $C$-edge $vw$ or $xy$, but not to any opposite minor 4- or 5-face \emph{(see Figure~\ref{fig:5faceD})}.}

Without loss of generality, let the middle $C$-edge of $f'$ be $xy$. Then $vy \notin E(G)$, as otherwise $C$ can be extended by replacing the path $(v,w,x,y,z)$ with $(v,y,x,w,d,z)$. Let $f_1$ be the face opposite to $f$ that has $vw$ as a $C$-edge. Since $w(f) <\nicefrac{10}{3}$, $f_1$ is a either a minor 3-face with middle $C$-edge $uv$ or a minor 2-face with $C$-edges $vw$ and $wx$. Assume to the contrary that $f_1$ is a 2-face. Then $vx \notin E(G)$, as otherwise $C$ can be extended by detouring through $b$. This implies $vz \in E(G)$, as otherwise $\deg_G(v) = 2$. Then $\{w,z\}$ is a 2-separator of $G$, which is a contradiction.

Hence, $f_1$ is a 3-face (see Figure~\ref{fig:5faceD}). Then $ux \notin E(G)$, as otherwise $C$ can be extended by replacing the path $(r,u,v,w,x)$ with $(r,b,w,v,u,x)$. Thus, since $\{w,z\}$ is no 2-separator of $G$, $uy$ or $vx$ is an edge of $G$. Assume to the contrary that $uy \notin E(G)$. Then $vx \in E(G)$, and we have $wy \notin E(G)$, as otherwise $C$ can be extended by replacing the path $(r,u,v,w,x,y,z)$ with $(r,b,w,y,x,v,u,a,z)$. Since $\deg_G(y) \geq 3$, this implies $uy \in E(G)$. Assume to the contrary that $vx \notin E(G)$. Then $xz \in E(G)$, as otherwise $\deg_G(x)=2$, and $C$ can be extended by replacing the path $(r,u,v,w,x,y,z)$ with $(r,b,w,v,u,y,x,z)$, which gives a contradiction. Hence, $uy \in E(G)$ and $vx \in E(G)$. Then $C$ can be extended by replacing the path $(u,v,w,x,y,z)$ with $(u,y,x,v,w,d,z)$.

\begin{figure}[!htb]
\begin{center}
	\begin{tikzpicture}[scale=0.25, -, 
	vertex/.style={circle,fill=black,draw,minimum size=5pt,inner sep=0pt}]
	\node[vertex,label=270:{ $r$}] (u) at (-12,0) [] {$ $};
	\node[vertex,label=270:{ $u$}] (v) at (-8,0) [] {$ $};
	\node[vertex,label=270:{ $v$}] (w) at (-4,0) [] {$ $};
	\node[vertex,label=270:{ $w$}] (x) at (0,0) [] {$ $};
	\node[vertex,label=270:{ $x$}] (y) at (4,0) [] {$ $};
	\node[vertex,label=270:{ $y$}] (z) at (8,0) [] {$ $};
	\node[vertex,label=270:{ $z$}] (q) at (12,0) [] {$ $};
	\node[vertex,label=90:{ $a$}] (a) at (2,7) [] {$ $};
	\node[vertex,label=-90:{ $b$}] (b) at (-6,-7) [] {$ $};
	\node[vertex,label=-90:{ $d$}] (d) at (6,-7) [] {$ $};
	\draw (2,4) node {$f$};
	\draw (6,-4) node {$f'$};
	\draw (-6,-4) node {$f_1$};
	\node (p1) at(-12,0) [] {$ $};
	\node (p2) at(12,0) [] {$ $};
	\node (q1) at(-16,1) [] {$ $};
	\node[label=0:{ $C$}] (q2) at(16,1) [] {$ $};
	\path[-,out=-30, in=180] (q1) edge (p1);
	\draw (u) -- (q);
	\path[-,out=0, in=-150] (p2) edge (q2);
	\draw (v) -- (a);
	\draw (q) -- (a);
	\path[-,out=0, in=-150] (p2) edge (q2);
	\path[dashed,out=35, in=145] (v) edge (z);
	\path[dashed,out=45, in=135] (w) edge (y);
	\draw (d) -- (x);
	\draw (d) -- (q);
	\draw (u) -- (b);
	\draw (b) -- (x);
	\end{tikzpicture}
\end{center}
\caption{Case~5d}\label{fig:5faceD}
\end{figure}

\item[Case 5e:] \textit{$f$ sends weight to an opposite minor 3-face $f'$ with middle $C$-edge $uv$ or $yz$, but not to any opposite minor 4- or 5-face \emph{(see Figure~\ref{fig:5faceE})}.}

Without loss of generality, let the middle $C$-edge of $f'$ be $yz$. Assume first that $f$ sends weight to a second opposite minor 3-face $f_1 \neq f'$. By Case~5d, $f_1$ has not middle $C$-edge $vw$, so that $f'$ must have middle $C$-edge $uv$. Then $wx$ is a $C$-edge of a major face opposite to $f$ that sends weight 1 to $f$, which contradicts $w(f) < \nicefrac{10}{3}$.

Hence, in order to satisfy $w(f) < \nicefrac{10}{3}$, $f$ sends by~R2 a total weight of $\nicefrac{4}{3}$ to opposite minor 2-faces. This implies that there is a minor 2-face $f_2$ opposite to $f$ that satisfies $m_{f,f_2}=2$. Then $f_2$ has not $C$-edges $uv$ and $vw$, as otherwise $wx$ would once again be a $C$-edge of a major face, which contradicts $w(f) < \nicefrac{10}{3}$. Hence, $f_2$ has $C$-edges $vw$ and $wx$ (see Figure~\ref{fig:5faceE}). Then $uw \notin E(G)$, as otherwise $C$ can be extended by replacing the path $(u,v,w,x)$ with $(u,w,v,b,x)$, and $wy \notin E(G)$, as otherwise $C$ can be extended by replacing the path $(v,w,x,y)$ with $(v,b,x,w,y)$. Since $\deg_G(w) \geq 3$, $wz \in E(G)$. Then $C$ can be extended by replacing the path $(w,x,y,z,q)$ with $(w,z,y,x,d,q)$, which is a contradiction.

\begin{figure}[!htb]
\begin{center}
	\begin{tikzpicture}[scale=0.25, -, 
	vertex/.style={circle,fill=black,draw,minimum size=5pt,inner sep=0pt}]
	\node[vertex,label=270:{ $u$}] (v) at (-8,0) [] {$ $};
	\node[vertex,label=270:{ $v$}] (w) at (-4,0) [] {$ $};
	\node[vertex,label=270:{ $w$}] (x) at (0,0) [] {$ $};
	\node[vertex,label=270:{ $x$}] (y) at (4,0) [] {$ $};
	\node[vertex,label=270:{ $y$}] (z) at (8,0) [] {$ $};
	\node[vertex,label=270:{ $z$}] (q) at (12,0) [] {$ $};
	\node[vertex,label=270:{ $q$}] (r) at (16,0) [] {$ $};
	\node[vertex,label=90:{ $a$}] (a) at (2,7) [] {$ $};
	\node[vertex,label=-90:{ $b$}] (b) at (0,-7) [] {$ $};
	\node[vertex,label=-90:{ $d$}] (d) at (10,-7) [] {$ $};
	\draw (2,4) node {$f$};
	\draw (0,-4) node {$f_2$};
	\draw (10,-4) node {$f'$};
	\node (p1) at(-8,0) [] {$ $};
	\node (p2) at(16,0) [] {$ $};
	\node (q1) at(-12,1) [] {$ $};
	\node[label=0:{ $C$}] (q2) at(20,1) [] {$ $};
	\path[-,out=-30, in=180] (q1) edge (p1);
	\draw (v) -- (r);
	\path[-,out=0, in=-150] (p2) edge (q2);
	\draw (v) -- (a);
	\draw (q) -- (a);
	\path[-,out=0, in=-150] (p2) edge (q2);
	\path[dashed,out=35, in=145] (x) edge (q);
	\draw (y) -- (b);
	\draw (b) -- (w);
	\draw (y) -- (d);
	\draw (d) -- (r);
	\end{tikzpicture}
\end{center}
\caption{Case~5e}\label{fig:5faceE}
\end{figure}
\end{description}

We conclude that $f$ sends no weight to any opposite minor 3-, 4- or 5-face. In order to satisfy $w(f) <\nicefrac{10}{3}$, $f$ must therefore send a total weight of \nicefrac{6}{3} to opposite minor 2-faces by~R2. In particular, there is at least one minor 2-face $f'$ opposite to $f$ that has $m_{f,f'}=2$. We distinguish the following subcases for $f'$.

\begin{description}
\item[Case 5f:] \textit{$f'$ has $C$-edges $uv$ and $vw$, or $xy$ and $yz$ \emph{(see Figure~\ref{fig:5faceF})}.}

Without loss of generality, let $f'$ have $C$-edges $xy$ and $yz$. Assume first that $f$ has a second opposite minor 2-face $f_1 \neq f'$ with $m_{f,f_1}=2$. Then $f_1$ has not $C$-edges $uv$ and $vw$, as then $wx$ would be a $C$-edge of a major face sending $f$ weight 1, which implies $w(f)=5-4 \cdot \frac23+1=\nicefrac{10}{3}$. Hence, $f_1$ has $C$-edges $vw$ and $wx$ (see Figure~\ref{fig:5faceF}). Then $wy \notin E(G)$, as otherwise $C$ can be extended by replacing the path $(w,x,y,z)$ with $(w,y,x,d,z)$. Hence, $vy \notin E(G)$, as otherwise $\deg_G(w) = 2$. Since $\deg_G(y) \geq 3$, we conclude $uy \in E(G)$ and, by $\deg_G(w) \geq 3$, $uw \in E(G)$. Then $C$ can be extended by replacing the path $(u,v,w,x)$ with $(u,w,v,b,x)$.

Hence, $f$ has no second opposite minor 2-face $f_1 \neq f'$ with $m_{f,f_1}=2$. Since $f$ sends a total weight of \nicefrac{6}{3} to opposite minor 2-faces by R2, $f$ has an opposite minor 2-face $f_2 \neq f'$ that has $C$-edge $uv$ but no other $C$-edge of $f$. Then $vw$ and $wx$ are $C$-edges of major face(s), which contradicts $w(f) <\nicefrac{10}{3}$.

\begin{figure}[!htb]
\begin{center}
	\begin{tikzpicture}[scale=0.25, -, 
	vertex/.style={circle,fill=black,draw,minimum size=5pt,inner sep=0pt}]
	\node[vertex,label=270:{ $u$}] (v) at (-8,0) [] {$ $};
	\node[vertex,label=270:{ $v$}] (w) at (-4,0) [] {$ $};
	\node[vertex,label=270:{ $w$}] (x) at (0,0) [] {$ $};
	\node[vertex,label=270:{ $x$}] (y) at (4,0) [] {$ $};
	\node[vertex,label=270:{ $y$}] (z) at (8,0) [] {$ $};
	\node[vertex,label=270:{ $z$}] (q) at (12,0) [] {$ $};
	\node[vertex,label=270:{ $q$}] (r) at (16,0) [] {$ $};
	\node[vertex,label=90:{ $a$}] (a) at (2,7) [] {$ $};
	\node[vertex,label=-90:{ $b$}] (b) at (0,-7) [] {$ $};
	\node[vertex,label=-90:{ $d$}] (d) at (8,-7) [] {$ $};
	\draw (2,4) node {$f$};
	\draw (0,-4) node {$f_1$};
	\draw (8,-4) node {$f'$};
	\node (p1) at(-8,0) [] {$ $};
	\node (p2) at(16,0) [] {$ $};
	\node (q1) at(-12,1) [] {$ $};
	\node[label=0:{ $C$}] (q2) at(20,1) [] {$ $};
	\path[-,out=-30, in=180] (q1) edge (p1);
	\draw (v) -- (r);
	\path[-,out=0, in=-150] (p2) edge (q2);
	\draw (v) -- (a);
	\draw (q) -- (a);
	\path[-,out=0, in=-150] (p2) edge (q2);
	\path[dashed,out=35, in=145] (v) edge (z);
	\path[dashed,out=35, in=145] (v) edge (x);
	\draw (y) -- (b);
	\draw (b) -- (w);
	\draw (y) -- (d);
	\draw (d) -- (q);
	\end{tikzpicture}
\end{center}
\caption{Case~5f}\label{fig:5faceF}
\end{figure}

\item[Case 5g:] \textit{$f'$ has $C$-edges $vw$ and $wx$, or $wx$ and $xy$ \emph{(see Figure~\ref{fig:5faceG})}.}

Without loss of generality, let $f'$ have $C$-edges $wx$ and $xy$. By Case~5f, $f$ has no second opposite minor 2-face $f_1 \neq f'$ with $m_{f,f_1}=2$. By $w(f)<\nicefrac{10}{3}$, $f$ has an opposite minor 2-face $f_2$ that has exactly one of the $C$-edges of $f$ as a $C$-edge. If this edge $e$ is not $yz$, $e = uv$ and then $vw$ is a $C$-edge of a major face, which contradicts $w(f) < \nicefrac{10}{3}$. Hence $e = yz$. Since neither $uv$ nor $vw$ is a $C$-edge of a major face, as this would again contradict $w(f) < \nicefrac{10}{3}$, $uv$ and $vw$ are $C$-edges of a minor $j$-face $f_3$ with $j \geq 4$ that does not receive any weight from $f$. Then $f_3$ sends weight \nicefrac{1}{3} to $f$ by~R5, which gives $w(f) = \nicefrac{10}{3}$ and thus a contradiction.

\begin{figure}[!htb]
\begin{center}
	\begin{tikzpicture}[scale=0.25, -, 
	vertex/.style={circle,fill=black,draw,minimum size=5pt,inner sep=0pt}]
	\node[vertex,label=270:{ $r$}] (u) at (-16,0) [] {$ $};
	\node[vertex,label=270:{ $s$}] (s) at (-12,0) [] {$ $};
	\node[vertex,label=270:{ $u$}] (v) at (-8,0) [] {$ $};
	\node[vertex,label=270:{ $v$}] (w) at (-4,0) [] {$ $};
	\node[vertex,label=270:{ $w$}] (x) at (0,0) [] {$ $};
	\node[vertex,label=270:{ $x$}] (y) at (4,0) [] {$ $};
	\node[vertex,label=270:{ $y$}] (z) at (8,0) [] {$ $};
	\node[vertex,label=270:{ $z$}] (q) at (12,0) [] {$ $};
	\node[vertex,label=270:{ $q$}] (r) at (16,0) [] {$ $};
	\node[vertex,label=90:{ $a$}] (a) at (2,7) [] {$ $};
	\node[vertex,label=-90:{ $b$}] (b) at (12,-7) [] {$ $};
	\node[vertex,label=-90:{ $d$}] (d) at (4,-7) [] {$ $};
	\node[vertex,label=-90:{ $g$}] (g) at (-8,-7) [] {$ $};
	\draw (2,4) node {$f$};
	\draw (12,-4) node {$f_2$};
	\draw (-8,-4) node {$f_3$};
	\draw (4,-4) node {$f'$};
	\node (p1) at(-16,0) [] {$ $};
	\node (p2) at(16,0) [] {$ $};
	\node (q1) at(-20,1) [] {$ $};
	\node[label=0:{ $C$}] (q2) at(20,1) [] {$ $};
	\path[-,out=-30, in=180] (q1) edge (p1);
	\draw (u) -- (r);
	\path[-,out=0, in=-150] (p2) edge (q2);
	\draw (v) -- (a);
	\draw (q) -- (a);
	\path[-,out=0, in=-150] (p2) edge (q2);
	\path[dashed,out=35, in=145] (v) edge (y);
	\path[dashed,out=-45, in=-135] (u) edge (w);
	\draw (z) -- (b);
	\draw (b) -- (r);
	\draw (x) -- (d);
	\draw (d) -- (z);
	\draw (g) -- (u);
	\draw (g) -- (x);
	\end{tikzpicture}
\end{center}
\caption{Case~5g}\label{fig:5faceG}
\end{figure}
\end{description}

\item[Case 6:] \textit{$f$ is a minor $6$-face \emph{(see Figure~\ref{fig:6face})}.}

Then $f$ is initially charged with weight 6. If $f$ looses a total net weight of at most \nicefrac{8}{3}, then $w(f) \geq \nicefrac{10}{3}$, so assume that $f$ looses a total net weight of at least $3$. We distinguish the following subcases.

\begin{figure}[!htb]
\begin{center}
	\begin{tikzpicture}[scale=0.25, -, 
	vertex/.style={circle,fill=black,draw,minimum size=5pt,inner sep=0pt}]
	\node[vertex,label=270:{ $t$}] (t) at (-16,0) [] {$ $};
	\node[vertex,label=270:{ $u$}] (u) at (-12,0) [] {$ $};
	\node[vertex,label=270:{ $v$}] (v) at (-8,0) [] {$ $};
	\node[vertex,label=270:{ $w$}] (w) at (-4,0) [] {$ $};
	\node[vertex,label=270:{ $x$}] (x) at (0,0) [] {$ $};
	\node[vertex,label=270:{ $y$}] (y) at (4,0) [] {$ $};
	\node[vertex,label=270:{ $z$}] (z) at (8,0) [] {$ $};
	\node[vertex,label=90:{ $a$}] (a) at (-4,7) [] {$ $};
	\draw (-4,4) node {$f$};
	\node (p1) at(-16,0) [] {$ $};
	\node (p2) at(8,0) [] {$ $};
	\node (q1) at(-20,1) [] {$ $};
	\node[label=0:{ $C$}] (q2) at(12,1) [] {$ $};
	\path[-,out=-30, in=180] (q1) edge (p1);
	\draw (t) -- (z);
	\path[-,out=0, in=-150] (p2) edge (q2);
	\draw (t) -- (a);
	\draw (z) -- (a);
	\end{tikzpicture}
\end{center}
\caption{Case~6}\label{fig:6face}
\end{figure}

\begin{description}
\item[Case 6a:] \textit{$f$ sends weight to an opposite minor 5-face $f'$ \emph{(see Figure~\ref{fig:6faceA})}.}

Without loss of generality, let $xy$ and $yz$ be $C$-edges of $f'$ getting weight from $f$ by~R5. Then $f$ sends weight \nicefrac{1}{3} to $f'$, and total weight \nicefrac{8}{3} to opposite minor 2-faces $f_3$ and $f_4$ by~R1--R5, as otherwise $w(f) \geq \nicefrac{10}{3}$ (see Figure~\ref{fig:6faceA}). Let $f_1$ and $f_2$ be the two minor 2-faces opposite to $f'$ due to~R5.

\begin{figure}[!htb]
\begin{center}
	\begin{tikzpicture}[scale=0.25, -, 
	vertex/.style={circle,fill=black,draw,minimum size=5pt,inner sep=0pt}]
	\node[vertex,label=270:{ $t$}] (t) at (-16,0) [] {$ $};
	\node[vertex,label=270:{ $u$}] (u) at (-12,0) [] {$ $};
	\node[vertex,label=270:{ $v$}] (v) at (-8,0) [] {$ $};
	\node[vertex,label=270:{ $w$}] (w) at (-4,0) [] {$ $};
	\node[vertex,label=270:{ $x$}] (x) at (0,0) [] {$ $};
	\node[vertex,label=270:{ $y$}] (y) at (4,0) [] {$ $};
	\node[vertex,label=270:{ $z$}] (z) at (8,0) [] {$ $};
	\node[vertex,label=270:{ $p$}] (p) at (12,0) [] {$ $};
	\node[vertex,label=270:{ $q$}] (q) at (16,0) [] {$ $};
	\node[vertex,label=270:{ $r$}] (r) at (20,0) [] {$ $};
	\node[vertex,label=270:{ $s$}] (s) at (24,0) [] {$ $};
	\node[vertex,label=90:{ $a$}] (a) at (-4,7) [] {$ $};
	\node[vertex,label=-90:{ $i$}] (i) at (-12,-7) [] {$ $};
	\node[vertex,label=-90:{ $b$}] (b) at (-4,-7) [] {$ $};
	\node[vertex,label=-90:{ $d$}] (d) at (10,-7) [] {$ $};
	\node[vertex,label=90:{ $g$}] (g) at (12,7) [] {$ $};
	\node[vertex,label=90:{ $h$}] (h) at (20,7) [] {$ $};
	\draw (-4,4) node {$f$};
	\draw (10,-4) node {$f'$};
	\draw (12,4) node {$f_1$};
	\draw (20,4) node {$f_2$};
	\draw (-12,-4) node {$f_3$};
	\draw (-4,-4) node {$f_4$};
	\node (p1) at(-16,0) [] {$ $};
	\node (p2) at(24,0) [] {$ $};
	\node (q1) at(-20,1) [] {$ $};
	\node[label=0:{ $C$}] (q2) at(28,1) [] {$ $};
	\path[-,out=-30, in=180] (q1) edge (p1);
	\draw (t) -- (s);
	\path[-,out=0, in=-150] (p2) edge (q2);
	\draw (t) -- (a);
	\draw (z) -- (a);
	\path[dashed,out=35, in=145] (w) edge (z);
	\path[dashed,out=-35, in=-145] (x) edge (p);
	\draw (x) -- (d);
	\draw (d) -- (r);
	\draw (v) -- (b);
	\draw (b) -- (x);
	\draw (z) -- (g);
	\draw (q) -- (g);
	\draw (q) -- (h);
	\draw (s) -- (h);
	\draw (i) -- (t);
	\draw (i) -- (v);
	\end{tikzpicture}
\end{center}
\caption{Case~6a}\label{fig:6faceA}
\end{figure}

We have $uw \notin E(G)$ and $wy \notin E(G)$, as otherwise $C$ can be extended by detouring through $b$, and $tw \notin E(G)$, as otherwise $\deg_G(u)=2$. Since $\deg_G(w) \geq 3$, $wz \in E(G)$. Moreover, $yp \notin E(G)$ and $pr \notin E(G)$, as otherwise $C$ can be extended by detouring through $g$. Since $\deg_G(p) \geq 3$, $xp \in E(G)$. Hence, $\deg_G(y) = 2$, which contradicts that $G$ is 3-connected.

\item[Case 6b:] \textit{$f$ sends weight to an opposite minor 4-face $f'$ \emph{(see Figure~\ref{fig:6faceB})}.}

Without loss of generality, let $xy$ and $yz$ be $C$-edges of $f'$ by~R4. Since $w(f)<\nicefrac{10}{3}$, $f$ has neither an opposite minor 5-face, nor a second opposite minor 4-face. Assume first that $f$ sends weight to an opposite minor 3-face $f_1$. Then $f$ sends total weight \nicefrac{5}{3} to $f'$ and $f_1$, and must therefore send weight \nicefrac{4}{3} to minor 2-face(s), as otherwise $w(f) \geq \nicefrac{10}{3}$. Hence, $f_1$ has middle $C$-edge $tu$, and $f$ has one opposite minor 2-face $f_2$ that has $C$-edges $vw$ and $wx$ (see Figure~\ref{fig:6faceB}).

\begin{figure}[!htb]
\begin{center}
	\begin{tikzpicture}[scale=0.25, -, 
	vertex/.style={circle,fill=black,draw,minimum size=5pt,inner sep=0pt}]
	\node[vertex,label=270:{ $s$}] (s) at (-20,0) [] {$ $};
	\node[vertex,label=270:{ $t$}] (t) at (-16,0) [] {$ $};
	\node[vertex,label=270:{ $u$}] (u) at (-12,0) [] {$ $};
	\node[vertex,label=270:{ $v$}] (v) at (-8,0) [] {$ $};
	\node[vertex,label=270:{ $w$}] (w) at (-4,0) [] {$ $};
	\node[vertex,label=270:{ $x$}] (x) at (0,0) [] {$ $};
	\node[vertex,label=270:{ $y$}] (y) at (4,0) [] {$ $};
	\node[vertex,label=270:{ $z$}] (z) at (8,0) [] {$ $};
	\node[vertex,label=270:{ $q$}] (q) at (12,0) [] {$ $};
	\node[vertex,label=270:{ $r$}] (r) at (16,0) [] {$ $};
	\node[vertex,label=90:{ $a$}] (a) at (-4,7) [] {$ $};
	\node[vertex,label=-90:{ $b$}] (b) at (-4,-7) [] {$ $};
	\node[vertex,label=-90:{ $d$}] (d) at (8,-7) [] {$ $};
	\node[vertex,label=90:{ $g$}] (g) at (12,7) [] {$ $};
	\node[vertex,label=-90:{ $h$}] (h) at (-14,-7) [] {$ $};
	\draw (-4,4) node {$f$};
	\draw (8,-4) node {$f'$};
	\draw (12,4) node {$f''$};
	\draw (-4,-4) node {$f_2$};
	\draw (-14,-4) node {$f_1$};
	\node (p1) at(-20,0) [] {$ $};
	\node (p2) at(16,0) [] {$ $};
	\node (q1) at(-24,1) [] {$ $};
	\node[label=0:{ $C$}] (q2) at(20,1) [] {$ $};
	\path[-,out=-30, in=180] (q1) edge (p1);
	\draw (s) -- (r);
	\path[-,out=0, in=-150] (p2) edge (q2);
	\draw (t) -- (a);
	\draw (z) -- (a);
	\path[-,out=0, in=-150] (p2) edge (q2);
	\path[dashed,out=35, in=145] (w) edge (z);
	\path[dashed,out=-45, in=-135] (x) edge (q);
	\draw (x) -- (d);
	\draw (d) -- (r);
	\draw (v) -- (b);
	\draw (b) -- (x);
	\draw (z) -- (g);
	\draw (g) -- (r);
	\draw (s) -- (h);
	\draw (v) -- (h);
	\end{tikzpicture}
\end{center}
\caption{Case~6b}\label{fig:6faceB}
\end{figure}

Then $uw$ and $wy$ are not edges of $G$, as otherwise $C$ can be extended by detouring through $b$. Moreover, $tw \notin E(G)$, as otherwise $C$ can be extended by replacing the path $(s,t,u,v,w)$ with $(s,h,v,u,t,w)$. Hence, $wz \in E(G)$, as otherwise $\deg_G(w)=2$.
Moreover, $yq \notin E(G)$ and $xq \in E(G)$ for the same reason as in Case~4i, which contradicts $\deg_G(y) \geq 3$.

\item[Case 6c:] \textit{$f$ sends weight to an opposite minor 3-face $f'$ with middle $C$-edge $vw$ or $wx$ \emph{(see Figure~\ref{fig:6faceC})}.}

Without loss of generality, let the middle $C$-edge of $f'$ be $wx$. In order to have $w(f) < \nicefrac{10}{3}$, $f$ must by~R2--R3 send weight $2$ to minor 2-faces. Thus, $f$ has two minor 2-faces $f_1$ and $f_2$ such that $f_1$ has $C$-edges $tu$ and $uv$, and $f_2$ has $yz$ as a $C$-edge.

\begin{figure}[!htb]
\begin{center}
	\begin{tikzpicture}[scale=0.25, -, 
	vertex/.style={circle,fill=black,draw,minimum size=5pt,inner sep=0pt}]
	\node[vertex,label=270:{ $t$}] (u) at (-12,0) [] {$ $};
	\node[vertex,label=270:{ $u$}] (v) at (-8,0) [] {$ $};
	\node[vertex,label=270:{ $v$}] (w) at (-4,0) [] {$ $};
	\node[vertex,label=270:{ $w$}] (x) at (0,0) [] {$ $};
	\node[vertex,label=270:{ $x$}] (y) at (4,0) [] {$ $};
	\node[vertex,label=270:{ $y$}] (z) at (8,0) [] {$ $};
	\node[vertex,label=270:{ $z$}] (q) at (12,0) [] {$ $};
	\node[vertex,label=270:{ $q$}] (r) at (16,0) [] {$ $};
	\node[vertex,label=90:{ $a$}] (a) at (0,7) [] {$ $};
	\node[vertex,label=-90:{ $b$}] (b) at (-8,-7) [] {$ $};
	\node[vertex,label=-90:{ $d$}] (d) at (2,-7) [] {$ $};
	\node[vertex,label=-90:{ $g$}] (g) at (12,-7) [] {$ $};
	\draw (0,4) node {$f$};
	\draw (2,-4) node {$f'$};
	\draw (-8,-4) node {$f_1$};
	\draw (12,-4) node {$f_2$};
	\node (p1) at(-12,0) [] {$ $};
	\node (p2) at(16,0) [] {$ $};
	\node (q1) at(-16,1) [] {$ $};
	\node[label=0:{ $C$}] (q2) at(20,1) [] {$ $};
	\path[-,out=-30, in=180] (q1) edge (p1);
	\draw (u) -- (r);
	\path[-,out=0, in=-150] (p2) edge (q2);
	\draw (u) -- (a);
	\draw (q) -- (a);
	\path[-,out=0, in=-150] (p2) edge (q2);
	\path[dashed,out=30, in=150] (v) edge (q);
	\draw (d) -- (z);
	\draw (d) -- (w);
	\draw (u) -- (b);
	\draw (b) -- (w);
	\draw (z) -- (g);
	\draw (g) -- (r);
	\end{tikzpicture}
\end{center}
\caption{Case~6c}\label{fig:6faceC}
\end{figure}

Then $uw \notin E(G)$, as otherwise $C$ can be extended by detouring $C$ through $b$. In addition, $ux \notin E(G)$, as otherwise $C$ can be extended by replacing the path $(u,v,w,x,y)$ with $(u,x,w,v,d,y)$. Then $uy \notin E(G)$, as otherwise the fact that $\{v,y\}$ is not a 2-separator of $G$ would imply $uw \in E(G)$ or $ux \in E(G)$. Since $\deg_G(u) \geq 3$, $uz \in E(G)$. Then we can obtain $\overline{C}$ from $C$ by replacing the path $(t,u,v,w,x,y,z,q)$ with $(t,a,z,u,v,w,x,y,g,q)$.

\item[Case 6d:] \textit{$f$ sends weight to an opposite minor 3-face $f'$ with middle $C$-edge $uv$ or $xy$ \emph{(see Figure~\ref{fig:6faceD})}.}

Without loss of generality, let the middle $C$-edge of $f'$ be $xy$. As in Case~6c, $w(f)<\nicefrac{10}{3}$ implies that $f$ has opposite minor 2-faces $f_1$ and $f_2$ such that $f_2$ has $C$-edges $uv$ and $vw$ and $f_1$ has $C$-edge $tu$ (see Figure~\ref{fig:6faceD}).

\begin{figure}[!htb]
\begin{center}
	\begin{tikzpicture}[scale=0.25, -, 
	vertex/.style={circle,fill=black,draw,minimum size=5pt,inner sep=0pt}]
	\node[vertex,label=270:{ $s$}] (s) at (-16,0) [] {$ $};
	\node[vertex,label=270:{ $t$}] (u) at (-12,0) [] {$ $};
	\node[vertex,label=270:{ $u$}] (v) at (-8,0) [] {$ $};
	\node[vertex,label=270:{ $v$}] (w) at (-4,0) [] {$ $};
	\node[vertex,label=270:{ $w$}] (x) at (0,0) [] {$ $};
	\node[vertex,label=270:{ $x$}] (y) at (4,0) [] {$ $};
	\node[vertex,label=270:{ $y$}] (z) at (8,0) [] {$ $};
	\node[vertex,label=270:{ $z$}] (q) at (12,0) [] {$ $};
	\node[vertex,label=90:{ $a$}] (a) at (0,7) [] {$ $};
	\node[vertex,label=-90:{ $b$}] (b) at (-4,-7) [] {$ $};
	\node[vertex,label=-90:{ $d$}] (d) at (6,-7) [] {$ $};
	\node[vertex,label=-90:{ $g$}] (g) at (-12,-7) [] {$ $};
	\draw (0,4) node {$f$};
	\draw (6,-4) node {$f'$};
	\draw (-12,-4) node {$f_1$};
	\draw (-4,-4) node {$f_2$};
	\node (p1) at(-16,0) [] {$ $};
	\node (p2) at(12,0) [] {$ $};
	\node (q1) at(-20,1) [] {$ $};
	\node[label=0:{ $C$}] (q2) at(16,1) [] {$ $};
	\path[-,out=-30, in=180] (q1) edge (p1);
	\draw (s) -- (q);
	\path[-,out=0, in=-150] (p2) edge (q2);
	\draw (u) -- (a);
	\draw (q) -- (a);
	\path[-,out=0, in=-150] (p2) edge (q2);
	\path[dashed,out=35, in=145] (w) edge (q);
	\draw (d) -- (x);
	\draw (d) -- (q);
	\draw (v) -- (b);
	\draw (b) -- (x);
	\draw (s) -- (g);
	\draw (v) -- (g);
	\end{tikzpicture}
\end{center}
\caption{Case~6d}\label{fig:6faceD}
\end{figure}

Then $tv$ and $vx$ are not edges of $G$, as otherwise $C$ can be extended by detouring $C$ through $b$. In addition, $vy \notin E(G)$, as otherwise $C$ can be extended by replacing the path $(v,w,x,y,z)$ with $(v,y,x,w,d,z)$. Since $\deg_G(v) \geq 3$, $vz \in E(G)$. This implies that $\{w,z\}$ is a 2-separator of $G$, which contradicts that $G$ is 3-connected.

\item[Case 6e:] \textit{$f$ sends weight to an opposite minor 3-face $f'$ with middle $C$-edge $tu$ or $yz$, but not to any opposite minor 4- or 5-face \emph{(see Figure~\ref{fig:6faceE})}.}

Without loss of generality, let the middle $C$-edge of $f'$ be $yz$. Assume first that $f$ has a second opposite minor 3-face $f''$. By Cases~6c+d, $f''$ has middle $C$-edge $tu$. By $w(f)<\nicefrac{10}{3}$, $f$ has an opposite minor 2-face $f_2$ with $C$-edges $vw$ and $wx$ (see Figure~\ref{fig:6faceE}). Then $uw \notin E(G)$ and $wy \notin E(G)$, as otherwise $C$ can be extended by detouring through $b$. Moreover, $wz \notin E(G)$, as otherwise $C$ can be extended by replacing the path $(w,x,y,z,q)$ with $(w,z,y,x,d,q)$. By symmetry, $tw \notin E(G)$, which contradicts $\deg_G(w) \geq 3$.

\begin{figure}[!htb]
\begin{center}
	\begin{tikzpicture}[scale=0.25, -, 
	vertex/.style={circle,fill=black,draw,minimum size=5pt,inner sep=0pt}]
	\node[vertex,label=270:{ $t$}] (u) at (-12,0) [] {$ $};
	\node[vertex,label=270:{ $u$}] (v) at (-8,0) [] {$ $};
	\node[vertex,label=270:{ $v$}] (w) at (-4,0) [] {$ $};
	\node[vertex,label=270:{ $w$}] (x) at (0,0) [] {$ $};
	\node[vertex,label=270:{ $x$}] (y) at (4,0) [] {$ $};
	\node[vertex,label=270:{ $y$}] (z) at (8,0) [] {$ $};
	\node[vertex,label=270:{ $z$}] (q) at (12,0) [] {$ $};
	\node[vertex,label=270:{ $q$}] (r) at (16,0) [] {$ $};
	\node[vertex,label=90:{ $a$}] (a) at (0,7) [] {$ $};
	\node[vertex,label=-90:{ $b$}] (b) at (0,-7) [] {$ $};
	\node[vertex,label=-90:{ $d$}] (d) at (10,-7) [] {$ $};
	\draw (0,4) node {$f$};
	\draw (0,-4) node {$f_2$};
	\draw (10,-4) node {$f'$};
	\node (p1) at(-12,0) [] {$ $};
	\node (p2) at(16,0) [] {$ $};
	\node (q1) at(-16,1) [] {$ $};
	\node[label=0:{ $C$}] (q2) at(20,1) [] {$ $};
	\path[-,out=-30, in=180] (q1) edge (p1);
	\draw (u) -- (r);
	\path[-,out=0, in=-150] (p2) edge (q2);
	\draw (u) -- (a);
	\draw (q) -- (a);
	\path[-,out=0, in=-150] (p2) edge (q2);
	\path[dashed,out=35, in=145] (u) edge (x);
	\draw (y) -- (b);
	\draw (b) -- (w);
	\draw (y) -- (d);
	\draw (d) -- (r);
	\end{tikzpicture}
\end{center}
\caption{Case~6e}\label{fig:6faceE}
\end{figure}

Hence, by~R1--R3, $f$ sends total weight $2$ to at least two opposite minor 2-faces $f_1$ and $f_2$. If $m_{f,f_1}=1$ or $m_{f,f_2}=1$, either the edge $uv$ or the edge $wx$ would be a $C$-edge of a major face, which contradicts $w(f) < \nicefrac{10}{3}$. Thus, $f_1$ has $C$-edges $tu$ and $uv$, and $f_2$ has $C$-edges $vw$ and $wx$. From the previous argument, we know that $uw$, $wy$ and $wz$ are not in $G$. Since $\deg_G(w) \geq 3$, $tw \in E(G)$. This contradicts $\deg_G(u) \geq 3$.
\end{description}

We conclude that $f$ sends no weight to any opposite minor 3-, 4- or 5-face. In order to satisfy $w(f) < \nicefrac{10}{3}$, $f$ must therefore send a total weight of \nicefrac{10}{3} to opposite minor 2-faces by~R2, as~R2 sends only multiples of weight \nicefrac{2}{3}. If some $C$-edge $e$ of $f$ is not a $C$-edge of a minor 2-face, $e$ must be either $tu$ or $yz$, as otherwise $e$ would be in a major face that sends weight 1 to $f$ and therefore contradicts $w(f)<\nicefrac{10}{3}$. Hence, $f$ has three opposite minor 2-faces $f_1$, $f_2$ and $f_3$ such that $m_{f,f_1} = m_{f,f_2} = 2$ and the $C$-edges of $f_1$ and $f_2$ are either $uv,vw,wx,xy$ or one of $tu,uv,vw,wx$ and $vw,wx,xy,yz$. We distinguish these subcases.

\begin{description}
\item[Case 6f:] \textit{The $C$-edges of $f_1$ and $f_2$ are $tu,uv,vw,wx$ or $vw,wx,xy,yz$ \emph{(see Figure~\ref{fig:6faceF})}.}

Without loss of generality, let $f_1$ and $f_2$ have the $C$-edges $vw,wx,xy,yz$. By the above argument, $f_3$ has the $C$-edges $tu$ and $uv$ (see Figure~\ref{fig:6faceF}).

\begin{figure}[!htb]
\begin{center}
	\begin{tikzpicture}[scale=0.25, -, 
	vertex/.style={circle,fill=black,draw,minimum size=5pt,inner sep=0pt}]
	\node[vertex,label=270:{ $t$}] (u) at (-12,0) [] {$ $};
	\node[vertex,label=270:{ $u$}] (v) at (-8,0) [] {$ $};
	\node[vertex,label=270:{ $v$}] (w) at (-4,0) [] {$ $};
	\node[vertex,label=270:{ $w$}] (x) at (0,0) [] {$ $};
	\node[vertex,label=270:{ $x$}] (y) at (4,0) [] {$ $};
	\node[vertex,label=270:{ $y$}] (z) at (8,0) [] {$ $};
	\node[vertex,label=270:{ $z$}] (q) at (12,0) [] {$ $};
	\node[vertex,label=90:{ $a$}] (a) at (0,7) [] {$ $};
	\node[vertex,label=-90:{ $b$}] (b) at (0,-7) [] {$ $};
	\node[vertex,label=-90:{ $d$}] (d) at (8,-7) [] {$ $};
	\node[vertex,label=-90:{ $g$}] (g) at (-8,-7) [] {$ $};
	\draw (0,4) node {$f$};
	\draw (-8,-4) node {$f_3$};
	\draw (0,-4) node {$f_1$};
	\draw (8,-4) node {$f_2$};
	\node (p1) at(-12,0) [] {$ $};
	\node (p2) at(12,0) [] {$ $};
	\node (q1) at(-16,1) [] {$ $};
	\node[label=0:{ $C$}] (q2) at(16,1) [] {$ $};
	\path[-,out=-30, in=180] (q1) edge (p1);
	\draw (u) -- (q);
	\path[-,out=0, in=-150] (p2) edge (q2);
	\draw (u) -- (a);
	\draw (q) -- (a);
	\path[-,out=0, in=-150] (p2) edge (q2);
	\draw (g) -- (u);
	\draw (g) -- (w);
	\draw (y) -- (b);
	\draw (b) -- (w);
	\draw (y) -- (d);
	\draw (d) -- (q);
	\end{tikzpicture}
\end{center}
\caption{Case~6f}\label{fig:6faceF}
\end{figure}

Then $uw$ and $wy$ are not in $G$, as otherwise $C$ can be extended by detouring through $b$. Moreoever, $wz \notin E(G)$, as otherwise $\deg_G(y) = 2$. By symmetry, $tw \notin E(G)$, which contradicts $\deg_G(w) \geq 3$.

\item[Case 6g:] \textit{The $C$-edges of $f_1$ and $f_2$ are $uv,vw,wx,xy$ \emph{(see Figure~\ref{fig:6faceG})}.}

Then $f_3$ has either $tu$ or $yz$ as a $C$-edge, say without loss of generality the latter.

\begin{figure}[!htb]
\begin{center}
	\begin{tikzpicture}[scale=0.25, -, 
	vertex/.style={circle,fill=black,draw,minimum size=5pt,inner sep=0pt}]
	\node[vertex,label=270:{ $t$}] (u) at (-12,0) [] {$ $};
	\node[vertex,label=270:{ $u$}] (v) at (-8,0) [] {$ $};
	\node[vertex,label=270:{ $v$}] (w) at (-4,0) [] {$ $};
	\node[vertex,label=270:{ $w$}] (x) at (0,0) [] {$ $};
	\node[vertex,label=270:{ $x$}] (y) at (4,0) [] {$ $};
	\node[vertex,label=270:{ $y$}] (z) at (8,0) [] {$ $};
	\node[vertex,label=270:{ $z$}] (q) at (12,0) [] {$ $};
	\node[vertex,label=270:{ $q$}] (r) at (16,0) [] {$ $};
	\node[vertex,label=90:{ $a$}] (a) at (0,7) [] {$ $};
	\node[vertex,label=-90:{ $b$}] (b) at (-4,-7) [] {$ $};
	\node[vertex,label=-90:{ $d$}] (d) at (4,-7) [] {$ $};
	\node[vertex,label=-90:{ $g$}] (g) at (12,-7) [] {$ $};
	\draw (0,4) node {$f$};
	\draw (12,-4) node {$f_3$};
	\draw (-4,-4) node {$f_1$};
	\draw (4,-4) node {$f_2$};
	\node (p1) at(-12,0) [] {$ $};
	\node (p2) at(16,0) [] {$ $};
	\node (q1) at(-16,1) [] {$ $};
	\node[label=0:{ $C$}] (q2) at(20,1) [] {$ $};
	\path[-,out=-30, in=180] (q1) edge (p1);
	\draw (u) -- (r);
	\path[-,out=0, in=-150] (p2) edge (q2);
	\draw (u) -- (a);
	\draw (q) -- (a);
	\path[-,out=0, in=-150] (p2) edge (q2);
	\path[dashed,out=35, in=145] (w) edge (q);
	\draw (g) -- (z);
	\draw (g) -- (r);
	\draw (v) -- (b);
	\draw (x) -- (b);
	\draw (x) -- (d);
	\draw (z) -- (d);
	\end{tikzpicture}
\end{center}
\caption{Case~6g}\label{fig:6faceG}
\end{figure}

Then $tv$ and $vx$ are not in $G$, as otherwise $C$ can be extended by detouring through $b$. Moreover, $vy \notin E(G)$, as otherwise $\deg_G(x) = 2$. Since $\deg_G(v) \geq 3$, $vz \in E(G)$. Then $xz \notin E(G)$, as otherwise $C$ can be extended by detouring through $g$. Hence, we obtain the contradiction $\deg_G(x) = 2$.
\end{description}

\item[Case 7:] \textit{$f$ is a minor $7$-face \emph{(see Figure~\ref{fig:7face})}.}

Then $f$ is initially charged with weight~7. If $f$ looses a total net weight of at most \nicefrac{11}{3}, then $w(f) \geq \nicefrac{10}{3}$, so assume that $f$ looses a total net weight of at least \nicefrac{12}{3}. According to R1--R5, $f$ sends to every opposite face $f'$ at most weight $\frac{2}{3}m_{f,f'}$ (for example, if $f'$ is a minor 3-face, $f$ sends only weight at most $\frac{1}{2}m_{f,f'}$ by~R3). Hence, $f$ does not send any weight to a 5-face, as otherwise $w(f) \geq \nicefrac{10}{3}$. We distinguish the remaining cases.

\begin{figure}[!htb]
\begin{center}
	\begin{tikzpicture}[scale=0.25, -, 
	vertex/.style={circle,fill=black,draw,minimum size=5pt,inner sep=0pt}]
	\node[vertex,label=270:{ $s$}] (s) at (-20,0) [] {$ $};
	\node[vertex,label=270:{ $t$}] (t) at (-16,0) [] {$ $};
	\node[vertex,label=270:{ $u$}] (u) at (-12,0) [] {$ $};
	\node[vertex,label=270:{ $v$}] (v) at (-8,0) [] {$ $};
	\node[vertex,label=270:{ $w$}] (w) at (-4,0) [] {$ $};
	\node[vertex,label=270:{ $x$}] (x) at (0,0) [] {$ $};
	\node[vertex,label=270:{ $y$}] (y) at (4,0) [] {$ $};
	\node[vertex,label=270:{ $z$}] (z) at (8,0) [] {$ $};
	\node[vertex,label=90:{ $a$}] (a) at (-6,7) [] {$ $};
	\draw (-6,4) node {$f$};
	\node (p1) at(-20,0) [] {$ $};
	\node (p2) at(8,0) [] {$ $};
	\node (q1) at(-24,1) [] {$ $};
	\node[label=0:{ $C$}] (q2) at(12,1) [] {$ $};
	\path[-,out=-30, in=180] (q1) edge (p1);
	\draw (s) -- (z);
	\path[-,out=0, in=-150] (p2) edge (q2);
	\draw (s) -- (a);
	\draw (z) -- (a);
	\end{tikzpicture}
\end{center}
\caption{Case~7}\label{fig:7face}
\end{figure}

\begin{description}
\item[Case 7a:] \textit{$f$ sends weight to an opposite minor 4-face $f'$ \emph{(see Figure~\ref{fig:7faceA})}.}

Without loss of generality, let $f'$ have $C$-edges $xy$ and $yz$. Since $w(f) < \nicefrac{10}{3}$, all other $C$-edges of $f$ are $C$-edges of minor 2-faces $f_1$, $f_2$ and $f_3$ (see Figure~\ref{fig:7faceA}).

\begin{figure}[!htb]
\begin{center}
	\begin{tikzpicture}[scale=0.25, -, 
	vertex/.style={circle,fill=black,draw,minimum size=5pt,inner sep=0pt}]
	\node[vertex,label=270:{ $r$}] (r) at (-24,0) [] {$ $};
	\node[vertex,label=270:{ $s$}] (s) at (-20,0) [] {$ $};
	\node[vertex,label=270:{ $t$}] (t) at (-16,0) [] {$ $};
	\node[vertex,label=270:{ $u$}] (u) at (-12,0) [] {$ $};
	\node[vertex,label=270:{ $v$}] (v) at (-8,0) [] {$ $};
	\node[vertex,label=270:{ $w$}] (w) at (-4,0) [] {$ $};
	\node[vertex,label=270:{ $x$}] (x) at (0,0) [] {$ $};
	\node[vertex,label=270:{ $y$}] (y) at (4,0) [] {$ $};
	\node[vertex,label=270:{ $z$}] (z) at (8,0) [] {$ $};
	\node[vertex,label=270:{ $p$}] (p) at (12,0) [] {$ $};
	\node[vertex,label=270:{ $q$}] (q) at (16,0) [] {$ $};
	\node[vertex,label=90:{ $a$}] (a) at (-4,7) [] {$ $};
	\node[vertex,label=-90:{ $b$}] (b) at (-4,-7) [] {$ $};
	\node[vertex,label=-90:{ $d$}] (d) at (8,-7) [] {$ $};
	\node[vertex,label=90:{ $g$}] (g) at (12,7) [] {$ $};
	\node[vertex,label=-90:{ $h$}] (h) at (-12,-7) [] {$ $};
	\node[vertex,label=-90:{ $i$}] (i) at (-20,-7) [] {$ $};
	\draw (-4,4) node {$f$};
	\draw (8,-4) node {$f'$};
	\draw (12,4) node {$f''$};
	\draw (-4,-4) node {$f_3$};
	\draw (-12,-4) node {$f_2$};
	\draw (-20,-4) node {$f_1$};
	\node (p1) at(-24,0) [] {$ $};
	\node (p2) at(16,0) [] {$ $};
	\node (q1) at(-28,1) [] {$ $};
	\node[label=0:{ $C$}] (q2) at(20,1) [] {$ $};
	\path[-,out=-30, in=180] (q1) edge (p1);
	\draw (r) -- (q);
	\path[-,out=0, in=-150] (p2) edge (q2);
	\draw (s) -- (a);
	\draw (z) -- (a);
	\path[-,out=0, in=-150] (p2) edge (q2);
	\path[dashed,out=25, in=155] (s) edge (w);
	\path[dashed,out=25, in=155] (s) edge (y);
	\path[dashed,out=-45, in=-135] (x) edge (p);
	\draw (x) -- (d);
	\draw (d) -- (q);
	\draw (v) -- (b);
	\draw (b) -- (x);
	\draw (z) -- (g);
	\draw (g) -- (q);
	\draw (t) -- (h);
	\draw (v) -- (h);
	\draw (r) -- (i);
	\draw (t) -- (i);
	\end{tikzpicture}
\end{center}
\caption{Case~7a}\label{fig:7faceA}
\end{figure}

Then $yp \notin E(G)$, as otherwise $C$ can be extended by detouring through $g$, and hence $xp \in E(G)$, as otherwise $\deg_G(p) = 2$. Also, $uw$ and $wy$ are not in $G$, as otherwise $C$ can be extended by detouring through $b$. Hence, $y$ has a neighbor in $G$ that is incident to $f$ and different from $\{w,x,z\}$. We conclude $wz \notin E(G)$. In addition, $tw \notin E(G)$, as otherwise $\deg_G(u) = 2$. Thus, $sw \in E(G)$, which implies $sy \in E(G)$. Then $\overline{C}$ can be obtained from $C$ by replacing the path $(r,s,t,u,v,w,x,y,z)$ with $(r,i,t,u,v,w,x,y,s,a,z)$.

\item[Case 7b:] \textit{$f$ sends weight to an opposite minor 3-face $f'$ \emph{(see Figure~\ref{fig:7faceB})}.}

Since $w(f)<\nicefrac{10}{3}$, the middle $C$-edge of $f'$ must be either $st$ or $yz$; say without loss of generality the latter. For the same reason as in Case~7a, all other $C$-edges of $f$ are $C$-edges of minor 2-faces $f_1$, $f_2$ and $f_3$ (see Figure~\ref{fig:7faceB}). 
Note that if there is another 3-face $f''$ with  middle $C$-edge $st$, then the edges $uv$, $vw$ and $wx$ are not all $C$-edges of some 2-face. 

\begin{figure}[!htb]
\begin{center}
	\begin{tikzpicture}[scale=0.25, -, 
	vertex/.style={circle,fill=black,draw,minimum size=5pt,inner sep=0pt}]
	\node[vertex,label=270:{ $r$}] (r) at (-20,0) [] {$ $};
	\node[vertex,label=270:{ $s$}] (s) at (-16,0) [] {$ $};
	\node[vertex,label=270:{ $t$}] (u) at (-12,0) [] {$ $};
	\node[vertex,label=270:{ $u$}] (v) at (-8,0) [] {$ $};
	\node[vertex,label=270:{ $v$}] (w) at (-4,0) [] {$ $};
	\node[vertex,label=270:{ $w$}] (x) at (0,0) [] {$ $};
	\node[vertex,label=270:{ $x$}] (y) at (4,0) [] {$ $};
	\node[vertex,label=270:{ $y$}] (z) at (8,0) [] {$ $};
	\node[vertex,label=270:{ $z$}] (q) at (12,0) [] {$ $};
	\node[vertex,label=270:{ $q$}] (p) at (16,0) [] {$ $};
	\node[vertex,label=90:{ $a$}] (a) at (-2,7) [] {$ $};
	\node[vertex,label=-90:{ $b$}] (b) at (0,-7) [] {$ $};
	\node[vertex,label=-90:{ $d$}] (d) at (10,-7) [] {$ $};
	\node[vertex,label=-90:{ $h$}] (h) at (-8,-7) [] {$ $};
	\node[vertex,label=-90:{ $g$}] (g) at (-16,-7) [] {$ $};
	\draw (-2,4) node {$f$};
	\draw (0,-4) node {$f_3$};
	\draw (-8,-4) node {$f_2$};
	\draw (-16,-4) node {$f_1$};
	\draw (10,-4) node {$f'$};
	\node (p1) at(-20,0) [] {$ $};
	\node (p2) at(16,0) [] {$ $};
	\node (q1) at(-24,1) [] {$ $};
	\node[label=0:{ $C$}] (q2) at(20,1) [] {$ $};
	\path[-,out=-30, in=180] (q1) edge (p1);
	\draw (r) -- (p);
	\path[-,out=0, in=-150] (p2) edge (q2);
	\draw (s) -- (a);
	\draw (q) -- (a);
	\path[-,out=0, in=-150] (p2) edge (q2);
	\path[dashed,out=30, in=150] (s) edge (x);
	\path[dashed,out=30, in=150] (s) edge (v);
	\draw (y) -- (b);
	\draw (b) -- (w);
	\draw (y) -- (d);
	\draw (d) -- (p);
	\draw (r) -- (g);
	\draw (u) -- (g);
	\draw (u) -- (h);
	\draw (w) -- (h);
	\end{tikzpicture}
\end{center}
\caption{Case~7b}\label{fig:7faceB}
\end{figure}

Then $uw \notin E(G)$ and $wy \notin E(G)$, as otherwise $C$ can be extended by detouring through $b$. Moreover, $wz \notin E(G)$, as otherwise $C$ can be extended by replacing the path $(w,x,y,z,q)$ with $(w,z,y,x,d,q)$. Also $tw \notin E(G)$, as otherwise $\deg_G(u) = 2$. Since $\deg_G(w) \geq 3$, $sw \in E(G)$. Since $\deg_G(u) \geq 3$, $su \in E(G)$. Then $C$ can be extended by replacing the path $(s,t,u,v)$ with $(s,u,t,h,v)$.

\item[Case 7c:] \textit{$f$ sends no weight to 3-, 4- and 5-faces \emph{(see Figure~\ref{fig:7faceC})}.}

Then $f$ sends a total weight of at least $6 \cdot \frac23 = 4$ to opposite minor 2-faces. The $C$-edges of these 2-faces must be consecutive on $C$, as otherwise exactly one $C$-edge of $f$ would be a $C$-edge of a major face, which contradicts $w(f)<\nicefrac{10}{3}$. Hence, there are three minor 2-faces $f_1$, $f_2$ and $f_3$, whose $C$-edges are consecutive on $C$ and satisfy $m_{f,f_1}=m_{f,f_2}=m_{f,f_3}=2$ (see Figure~\ref{fig:7faceC}). Assume without loss of generality that $f_3$ has $C$-edges $xy$ and $yz$.

\begin{figure}[!htb]
\begin{center}
	\begin{tikzpicture}[scale=0.25, -, 
	vertex/.style={circle,fill=black,draw,minimum size=5pt,inner sep=0pt}]
	\node[vertex,label=270:{ $s$}] (u) at (-12,0) [] {$ $};
	\node[vertex,label=270:{ $t$}] (v) at (-8,0) [] {$ $};
	\node[vertex,label=270:{ $u$}] (w) at (-4,0) [] {$ $};
	\node[vertex,label=270:{ $v$}] (x) at (0,0) [] {$ $};
	\node[vertex,label=270:{ $w$}] (y) at (4,0) [] {$ $};
	\node[vertex,label=270:{ $x$}] (z) at (8,0) [] {$ $};
	\node[vertex,label=270:{ $y$}] (q) at (12,0) [] {$ $};
	\node[vertex,label=270:{ $z$}] (r) at (16,0) [] {$ $};
	\node[vertex,label=90:{ $a$}] (a) at (2,7) [] {$ $};
	\node[vertex,label=-90:{ $b$}] (b) at (-4,-7) [] {$ $};
	\node[vertex,label=-90:{ $d$}] (d) at (4,-7) [] {$ $};
	\node[vertex,label=-90:{ $g$}] (g) at (12,-7) [] {$ $};
	\draw (2,4) node {$f$};
	\draw (12,-4) node {$f_3$};
	\draw (-4,-4) node {$f_1$};
	\draw (4,-4) node {$f_2$};
	\node (p1) at(-12,0) [] {$ $};
	\node (p2) at(16,0) [] {$ $};
	\node (q1) at(-16,1) [] {$ $};
	\node[label=0:{ $C$}] (q2) at(20,1) [] {$ $};
	\path[-,out=-30, in=180] (q1) edge (p1);
	\draw (u) -- (r);
	\path[-,out=0, in=-150] (p2) edge (q2);
	\draw (u) -- (a);
	\draw (r) -- (a);
	\path[-,out=0, in=-150] (p2) edge (q2);
	\path[dashed,out=30, in=150] (u) edge (y);
	\draw (g) -- (z);
	\draw (g) -- (r);
	\draw (v) -- (b);
	\draw (x) -- (b);
	\draw (x) -- (d);
	\draw (z) -- (d);
	\end{tikzpicture}
\end{center}
\caption{Case~7c}\label{fig:7faceC}
\end{figure}

Then $uw$ and $wy$ are not in $G$, as otherwise $C$ can be extended by detouring through $d$. Moreover, $tw$ and $wz$ are not in $G$, as otherwise $\deg_G(u) = 2$ or $\deg_G(y) = 2$. Since $\deg_G(w) \geq 3$, $sw \in E(G)$. Moreover, $su \notin E(G)$, as otherwise $C$ can be extended by detouring through $b$. Hence, we obtain the contradiction $\deg_G(u) = 2$.
\end{description}

\item[Case 8:] \textit{$f$ is a minor $8$-face \emph{(see Figure~\ref{fig:8face})}.}

Then $f$ is initially charged with weight~8. If $f$ looses a total net weight of at most \nicefrac{14}{3}, then $w(f) \geq \nicefrac{10}{3}$, so assume that $f$ looses a total net weight of at least \nicefrac{15}{3}. Hence, $f$ does not send any weight to a 4- or~5-face, as otherwise $w(f) \geq \nicefrac{10}{3}$. We distinguish the remaining cases.

\begin{figure}[!htb]
\begin{center}
	\begin{tikzpicture}[scale=0.25, -, 
	vertex/.style={circle,fill=black,draw,minimum size=5pt,inner sep=0pt}]
	\node[vertex,label=270:{ $r$}] (r) at (-24,0) [] {$ $};
	\node[vertex,label=270:{ $s$}] (s) at (-20,0) [] {$ $};
	\node[vertex,label=270:{ $t$}] (t) at (-16,0) [] {$ $};
	\node[vertex,label=270:{ $u$}] (u) at (-12,0) [] {$ $};
	\node[vertex,label=270:{ $v$}] (v) at (-8,0) [] {$ $};
	\node[vertex,label=270:{ $w$}] (w) at (-4,0) [] {$ $};
	\node[vertex,label=270:{ $x$}] (x) at (0,0) [] {$ $};
	\node[vertex,label=270:{ $y$}] (y) at (4,0) [] {$ $};
	\node[vertex,label=270:{ $z$}] (z) at (8,0) [] {$ $};
	\node[vertex,label=90:{ $a$}] (a) at (-8,7) [] {$ $};
	\draw (-8,4) node {$f$};
	\node (p1) at(-24,0) [] {$ $};
	\node (p2) at(8,0) [] {$ $};
	\node (q1) at(-28,1) [] {$ $};
	\node[label=0:{ $C$}] (q2) at(12,1) [] {$ $};
	\path[-,out=-30, in=180] (q1) edge (p1);
	\draw (r) -- (z);
	\path[-,out=0, in=-150] (p2) edge (q2);
	\draw (r) -- (a);
	\draw (z) -- (a);
	\end{tikzpicture}
\end{center}
\caption{Case~8}\label{fig:8face}
\end{figure}

\begin{description}
\item[Case 8a:] \textit{$f$ sends weight to an opposite minor 3-face $f'$ \emph{(see Figure~\ref{fig:8faceA})}.}

Then $w(f) < \nicefrac{10}{3}$ implies that $f'$ has exactly two $C$-edges that are $C$-edges of $f$, and that every other $C$-edge of $f$ is a $C$-edge of a minor 2-face. Without loss of generality, let $f'$ have middle $C$-edge $yz$, and let $f_1$, $f_2$ and $f_3$ be the minor 2-faces opposite to $f$ (see Figure~\ref{fig:8faceA}).

\begin{figure}[!htb]
\begin{center}
	\begin{tikzpicture}[scale=0.25, -, 
	vertex/.style={circle,fill=black,draw,minimum size=5pt,inner sep=0pt}]
	\node[vertex,label=270:{ $r$}] (r) at (-20,0) [] {$ $};
	\node[vertex,label=270:{ $s$}] (s) at (-16,0) [] {$ $};
	\node[vertex,label=270:{ $t$}] (u) at (-12,0) [] {$ $};
	\node[vertex,label=270:{ $u$}] (v) at (-8,0) [] {$ $};
	\node[vertex,label=270:{ $v$}] (w) at (-4,0) [] {$ $};
	\node[vertex,label=270:{ $w$}] (x) at (0,0) [] {$ $};
	\node[vertex,label=270:{ $x$}] (y) at (4,0) [] {$ $};
	\node[vertex,label=270:{ $y$}] (z) at (8,0) [] {$ $};
	\node[vertex,label=270:{ $z$}] (q) at (12,0) [] {$ $};
	\node[vertex,label=270:{ $q$}] (p) at (16,0) [] {$ $};
	\node[vertex,label=90:{ $a$}] (a) at (-4,7) [] {$ $};
	\node[vertex,label=-90:{ $b$}] (b) at (0,-7) [] {$ $};
	\node[vertex,label=-90:{ $d$}] (d) at (10,-7) [] {$ $};
	\node[vertex,label=-90:{ $h$}] (h) at (-8,-7) [] {$ $};
	\node[vertex,label=-90:{ $g$}] (g) at (-16,-7) [] {$ $};
	\draw (-4,4) node {$f$};
	\draw (0,-4) node {$f_3$};
	\draw (-8,-4) node {$f_2$};
	\draw (-16,-4) node {$f_1$};
	\draw (10,-4) node {$f'$};
	\node (p1) at(-20,0) [] {$ $};
	\node (p2) at(16,0) [] {$ $};
	\node (q1) at(-24,1) [] {$ $};
	\node[label=0:{ $C$}] (q2) at(20,1) [] {$ $};
	\path[-,out=-30, in=180] (q1) edge (p1);
	\draw (r) -- (p);
	\path[-,out=0, in=-150] (p2) edge (q2);
	\draw (r) -- (a);
	\draw (q) -- (a);
	\path[-,out=0, in=-150] (p2) edge (q2);
	\path[dashed,out=25, in=155] (r) edge (x);
	\path[dashed,out=25, in=155] (r) edge (v);
	\draw (y) -- (b);
	\draw (b) -- (w);
	\draw (y) -- (d);
	\draw (d) -- (p);
	\draw (r) -- (g);
	\draw (u) -- (g);
	\draw (u) -- (h);
	\draw (w) -- (h);
	\end{tikzpicture}
\end{center}
\caption{Case~8a}\label{fig:8faceA}
\end{figure}

Then $su$, $uw$ and $wy$ are not edges of $G$, as otherwise $C$ can be extended by detouring through $h$ or $b$. Moreover, $wz \notin E(G)$, as otherwise $C$ can be extended by replacing the path $(w,x,y,z,q)$ with $(w,z,y,x,d,q)$. Also $sw \notin E(G)$ and $tw \notin E(G)$, as otherwise $\deg_G(u) = 2$. Since $\deg_G(w) \geq 3$, $rw \in E(G)$. Since $\deg_G(u) \geq 3$, $ru \in E(G)$. This gives the contradiction $\deg_G(s) = 2$.

\item[Case 8b:] \textit{$f$ sends no weight to 3-, 4- and 5-faces \emph{(see Figure~\ref{fig:8faceB})}.}

Then $f$ sends a total weight of exactly $8 \cdot \frac23 = \nicefrac{16}{3}$ to opposite minor 2-faces, as~R2 sends only multiples of weight $\frac23$. Assume first that an opposite minor 2-face $f_4$ to $f$ has $C$-edges $xy$ and $yz$ (see Figure~\ref{fig:8faceB}). Then $wy \notin E(G)$, as otherwise $C$ can be extended by detouring through $g$, and $wz \notin E(G)$, as otherwise $\deg_G(y) = 2$. Then the same arguments as in Case~8a give the contradiction $\deg_G(s) = 2$.

\begin{figure}[!htb]
\begin{center}
	\begin{tikzpicture}[scale=0.25, -, 
	vertex/.style={circle,fill=black,draw,minimum size=5pt,inner sep=0pt}]
	\node[vertex,label=270:{ $r$}] (p) at (-16,0) [] {$ $};
	\node[vertex,label=270:{ $s$}] (u) at (-12,0) [] {$ $};
	\node[vertex,label=270:{ $t$}] (v) at (-8,0) [] {$ $};
	\node[vertex,label=270:{ $u$}] (w) at (-4,0) [] {$ $};
	\node[vertex,label=270:{ $v$}] (x) at (0,0) [] {$ $};
	\node[vertex,label=270:{ $w$}] (y) at (4,0) [] {$ $};
	\node[vertex,label=270:{ $x$}] (z) at (8,0) [] {$ $};
	\node[vertex,label=270:{ $y$}] (q) at (12,0) [] {$ $};
	\node[vertex,label=270:{ $z$}] (r) at (16,0) [] {$ $};
	\node[vertex,label=90:{ $a$}] (a) at (0,7) [] {$ $};
	\node[vertex,label=-90:{ $b$}] (b) at (-12,-7) [] {$ $};
	\node[vertex,label=-90:{ $d$}] (d) at (-4,-7) [] {$ $};
	\node[vertex,label=-90:{ $g$}] (g) at (4,-7) [] {$ $};
	\node[vertex,label=-90:{ $h$}] (h) at (12,-7) [] {$ $};
	\draw (0,4) node {$f$};
	\draw (-12,-4) node {$f_1$};
	\draw (-4,-4) node {$f_2$};
	\draw (4,-4) node {$f_3$};
	\draw (12,-4) node {$f_4$};
	\node (p1) at(-16,0) [] {$ $};
	\node (p2) at(16,0) [] {$ $};
	\node (q1) at(-20,1) [] {$ $};
	\node[label=0:{ $C$}] (q2) at(20,1) [] {$ $};
	\path[-,out=-30, in=180] (q1) edge (p1);
	\draw (p) -- (r);
	\path[-,out=0, in=-150] (p2) edge (q2);
	\draw (p) -- (a);
	\draw (r) -- (a);
	\path[-,out=0, in=-150] (p2) edge (q2);
	\path[dashed,out=25, in=155] (p) edge (w);
	\path[dashed,out=25, in=155] (p) edge (y);
	\draw (b) -- (p);
	\draw (b) -- (v);
	\draw (h) -- (z);
	\draw (h) -- (r);
	\draw (v) -- (d);
	\draw (x) -- (d);
	\draw (x) -- (g);
	\draw (z) -- (g);
	\end{tikzpicture}
\end{center}
\caption{Case~8b}\label{fig:8faceB}
\end{figure}

Hence, let $yz$ be the only $C$-edge of $f_4$ that is a $C$-edge of $f$. Then $v$ has no neighbor that is incident to $f$ and not in $\{u,w\}$, as otherwise $t$ or $x$ has degree~2 in $G$. Hence, we obtain the contradiction $\deg_G(v) = 2$.
\end{description}

\item[Case 9:] \textit{$f$ is a minor $j$-face with $j \geq 9$ \emph{(see Figure~\ref{fig:9faceA})}.}

Then $f$ is initially charged with weight $j$ and looses a total net weight of at most $\frac{2}{3}j$, so that $w(f) \geq \frac{1}{3}j \geq \frac{10}{3}$ if $j \geq 10$. Hence, $j=9$ and every $C$-edge of $f$ is a $C$-edge of a minor 2-face. Since 9 is odd, we may assume without loss of generality that one minor 2-face $f_1$ has $qr$ but no other $C$-edge of $f$ as a $C$-edge (see Figure~\ref{fig:9faceA}). Then the same arguments as in Cases~8a+b imply that $\deg_G(s) = 2$.

\begin{figure}[!htb]
\begin{center}
	\begin{tikzpicture}[scale=0.25, -, 
	vertex/.style={circle,fill=black,draw,minimum size=5pt,inner sep=0pt}]
	\node[vertex,label=270:{ $l$}] (l) at (-24,0) [] {$ $};
	\node[vertex,label=270:{ $q$}] (o) at (-20,0) [] {$ $};
	\node[vertex,label=270:{ $r$}] (p) at (-16,0) [] {$ $};
	\node[vertex,label=270:{ $s$}] (u) at (-12,0) [] {$ $};
	\node[vertex,label=270:{ $t$}] (v) at (-8,0) [] {$ $};
	\node[vertex,label=270:{ $u$}] (w) at (-4,0) [] {$ $};
	\node[vertex,label=270:{ $v$}] (x) at (0,0) [] {$ $};
	\node[vertex,label=270:{ $w$}] (y) at (4,0) [] {$ $};
	\node[vertex,label=270:{ $x$}] (z) at (8,0) [] {$ $};
	\node[vertex,label=270:{ $y$}] (q) at (12,0) [] {$ $};
	\node[vertex,label=270:{ $z$}] (r) at (16,0) [] {$ $};
	\node[vertex,label=90:{ $a$}] (a) at (0,7) [] {$ $};
	\node[vertex,label=-90:{ $b$}] (b) at (-12,-7) [] {$ $};
	\node[vertex,label=-90:{ $d$}] (d) at (-4,-7) [] {$ $};
	\node[vertex,label=-90:{ $g$}] (g) at (4,-7) [] {$ $};
	\node[vertex,label=-90:{ $h$}] (h) at (12,-7) [] {$ $};
	\node[vertex,label=-90:{ $i$}] (i) at (-20,-7) [] {$ $};
	\draw (0,4) node {$f$};
	\draw (-20,-4) node {$f_1$};
	\draw (-12,-4) node {$f_2$};
	\draw (-4,-4) node {$f_3$};
	\draw (4,-4) node {$f_4$};
	\draw (12,-4) node {$f_5$};
	\node (p1) at(-24,0) [] {$ $};
	\node (p2) at(16,0) [] {$ $};
	\node (q1) at(-28,1) [] {$ $};
	\node[label=0:{ $C$}] (q2) at(20,1) [] {$ $};
	\path[-,out=-30, in=180] (q1) edge (p1);
	\draw (l) -- (r);
	\path[-,out=0, in=-150] (p2) edge (q2);
	\draw (o) -- (a);
	\draw (r) -- (a);
	\path[-,out=0, in=-150] (p2) edge (q2);
	\path[dashed,out=20, in=160] (o) edge (w);
	\path[dashed,out=20, in=160] (o) edge (y);
	\draw (i) -- (l);
	\draw (i) -- (p);
	\draw (b) -- (p);
	\draw (b) -- (v);
	\draw (h) -- (z);
	\draw (h) -- (r);
	\draw (v) -- (d);
	\draw (x) -- (d);
	\draw (x) -- (g);
	\draw (z) -- (g);
	\end{tikzpicture}
\end{center}
\caption{Case~9}\label{fig:9faceA}
\end{figure}
\end{description}
This proves $2c=\sum_{f\in F(H)} w(f) \geq \nicefrac{10}{3}\cdot |M^- \cup M^+|$, which completes the proof of Theorem~\ref{thm:main}.
\hfill$\square$\\

\section{Remarks}\label{sec:remarks}
We remark that the bound of Theorem~\ref{thm:main} can be improved to $\frac{5}{8}(n+4)$ for every $n \geq 16$: then Lemma~5 in~\cite{Fabrici2016} implies the improved bound for the special case that $V^-$ or $V^+$ is empty, while in the remaining case $|V^-| \geq 1 \leq |V^+|$ Lemma~\ref{lem:minor} can be immediately strengthened to $|M^- \cup M^+| \geq |V^- \cup V^+| + 4$ using the same proof with a different induction base (see also~\cite{Fabrici2018}). This immediately improves the bound $\circ(G)\ge \frac{13}{21}(n+4)$ given in~\cite{Fabrici2016} for every $n \geq 16$. We note that $\circ(G) \geq \frac{5}{8}(n+4)$ does not hold for $n \leq 6$, as for these values a cycle of length at least $\frac{5}{8}(n+4) > n$ is impossible.

The proof of Theorem~\ref{thm:main} is constructive and gives a quadratic-time algorithm that finds a cycle of length at least $\frac{5}{8}(n+2)$, by applying the result of~\cite{Schmid2018b} exactly as shown in~\cite[Section Algorithm]{Fabrici2018}. We therefore conclude the following theorem.

\begin{theorem}
For every essentially 4-connected plane graph $G$ on $n$ vertices, a cycle of length at least $\frac{5}{8}(n+2)$ can be computed in time $O(n^2)$.
\end{theorem}

\bibliographystyle{abbrv}
\bibliography{Paper}

\begin{thebibliography}{1}

\bibitem{Dillencourt1996}
M.~B. Dillencourt.
\newblock Polyhedra of small order and their {H}amiltonian properties.
\newblock {\em Journal of Combinatorial Theory, Series B}, 66(1):87--122, 1996.

\bibitem{Fabrici2016}
I.~Fabrici, J.~Harant, and S.~Jendrol.
\newblock On longest cycles in essentially 4-connected planar graphs.
\newblock {\em Discussiones Mathematicae Graph Theory}, 36:565--575, 2016.

\bibitem{Fabrici2018}
I.~Fabrici, J.~Harant, S.~Mohr, and J.~M. Schmidt.
\newblock Longer cycles in essentially 4-connected planar graphs.
\newblock {\em Discussiones Mathematicae Graph Theory}, in press.

\bibitem{Grunbaum1976}
B.~Gr\"unbaum and J.~Malkevitch.
\newblock Pairs of edge-disjoint {H}amilton circuits.
\newblock {\em Aequationes Mathematicae}, 14:191--196, 1976.

\bibitem{Jackson1992}
B.~Jackson and N.~C. Wormald.
\newblock Longest cycles in 3-connected planar graphs.
\newblock {\em Journal of Combinatorial Theory, Series B}, 54:291--321, 1992.

\bibitem{Schmid2018b}
A.~Schmid and J.~M. Schmidt.
\newblock Computing {T}utte paths.
\newblock In {\em Proceedings of the 45th International Colloquium on Automata,
  Languages and Programming (ICALP'18)}, pages 98:1--98:14, 2018.

\bibitem{Whitney1932}
H.~Whitney.
\newblock Congruent graphs and the connectivity of graphs.
\newblock {\em American Journal of Mathematics}, 54(1):150--168, 1932.

\bibitem{Zhang1987}
C.-Q. Zhang.
\newblock Longest cycles and their chords.
\newblock {\em Journal of Graph Theory}, 11:341--345, 1987.

\end{thebibliography}
\end{document}